\documentclass[aap,MSNbibl,seceqn,dvips]{arximspdf}
\usepackage{url,breakurl}

\doi{10.1214/13-AAP941} 
\volume{24}
\issue{3}
\pubyear{2014}
\firstpage{1049}
\lastpage{1080}

\makeatletter
\newcommand{\rrVert}{\Vert}
\newcommand{\rrvert}{\vert}
\newcommand{\llVert}{\Vert}
\newcommand{\llvert}{\vert}
\newtheorem{theorem}{Theorem}[section]
\newtheorem{alem}[theorem]{Lemma}
\newtheorem{aprop}[theorem]{Proposition}
\newtheorem{acor}[theorem]{Corollary}
\newproclaim{arem}[theorem]{Remark}
\newtheorem{ahyp}[theorem]{Hypothesis}
\makeatother

\begin{document}
\begin{frontmatter}

\title{Pathwise optimal transport bounds between a~one-dimensional
diffusion and its Euler scheme}
\runtitle{Pathwise optimal transport bounds}

\begin{aug}
\author[A]{\fnms{A.} \snm{Alfonsi}\corref{}\ead[label=e1]{alfonsi@cermics.enpc.fr}\thanksref{t1}},
\author[A]{\fnms{B.} \snm{Jourdain}\ead[label=e2]{jourdain@cermics.enpc.fr}\thanksref{t1}}
\and
\author[B]{\fnms{A.} \snm{Kohatsu-Higa}\ead[label=e3]{arturokohatsu@gmail.com}\thanksref{t2}}
\runauthor{A. Alfonsi, B. Jourdain and A. Kohatsu-Higa}
\affiliation{Universit\'e Paris-Est, Universit\'e Paris-Est, and
Ritsumeikan University and~Japan Science and Technology Agency}
\address[A]{A. Alfonsi\\
B. Jourdain\\
CERMICS\\
Projet MathFi ENPC-INRIA-UMLV\\
Universit\'e Paris-Est\\
6-8 Avenue Blaise Pascal\\
77455 Marne La Vall\'ee, Cedex 2\\
France\\
\printead{e1}\\
\phantom{E-mail: }\printead*{e2}}
\address[B]{A. Kohatsu-Higa\\
Department of Mathematical Sciences\\
Ritsumeikan University\\
1-1-1 Nojihigashi\\
Kusatsu, Shiga 525-8577\\
Japan\\
and\\
Japan Science and Technology Agency\\
\printead{e3}}
\thankstext{t1}{Benefited from the support of the ``Chaire Risques
Financiers,'' Fondation du
Risque and of the French National Research Agency (ANR) under the program
ANR-08-BLAN-0218 BigMC and the Labex Bezout.}
\thankstext{t2}{Supported by grants of the Japanese goverment.}
\end{aug}

\received{\smonth{9} \syear{2012}}
\revised{\smonth{5} \syear{2013}}

%
\begin{abstract}
In the present paper, we prove that the Wasserstein distance on the
space of continuous sample-paths equipped with the supremum norm
between the laws of a uniformly elliptic one-dimensional diffusion
process and its Euler discretization with $N$ steps is smaller than
$O(N^{-2/3+\varepsilon})$ where $\varepsilon$ is an arbitrary positive
constant. This rate is intermediate between the strong error estimation
in $O(N^{-1/2})$ obtained when coupling the stochastic differential
equation and the Euler scheme with the same Brownian motion and the
weak error estimation $O(N^{-1})$ obtained when comparing the
expectations of the same function of the diffusion and of the Euler
scheme at the terminal time $T$. We also check that the supremum over
$t\in[0,T]$ of the Wasserstein distance on the space of probability
measures on the real line between the laws of the diffusion at time $t$
and the Euler scheme at time $t$ behaves like
$O(\sqrt{\log(N)}N^{-1})$.
\end{abstract}

%
\begin{keyword}[class=AMS]
\kwd{65C30}
\kwd{60H35}
\end{keyword}
\begin{keyword}
\kwd{Euler scheme}
\kwd{Wasserstein distance}
\kwd{weak trajectorial error}
\kwd{diffusion bridges}
\end{keyword}

\end{frontmatter}

\setcounter{section}{-1}
\section{Introduction}

For $\sigma\dvtx \mathbb{R}\to\mathbb{R}$ and $b\dvtx
\mathbb{R}\to\mathbb {R}$, we are interested in the simulation of the
stochastic differential equation
%
%
\begin{equation}
dX_t=\sigma(X_t)\,dW_t+b(X_t)\,dt,
\label{sde}
\end{equation}
where $X_0=x_0\in\mathbb{R}$ and $W=(W_t)_{t\geq0}$ is a standard Brownian
motion. We make the standard Lipschitz assumptions on the coefficients,
\[
\exists K\in(0,+\infty), \forall x,y\in\mathbb{R}\qquad
\bigl|\sigma(x)-\sigma(y)\bigr|+\bigl|b(x)-b(y)\bigr|\leq K|x-y|.
\]
%


For $T>0$, we are interested in the approximation of
$X=(X_t)_{t\in[0,T]}$ by its Euler scheme
$\bar{X}=(\bar{X}_t)_{t\in[0,T]}$ with $N\geq1$ time-steps. We consider
the regular grid $\{0=t_0<t_1<t_2<\cdots<t_N=T\}$ of the interval
$[0,T]$ with $t_k=\frac{kT}{N}$ and define inductively $\bar{X}_0=x_0$
and
%
%
\begin{equation}
\bar{X}_t=\bar{X}_{t_k}+\sigma(\bar{X}_{t_k})
(W_t-W_{t_k})+b(\bar{X}_{t_k})
(t-t_k)\qquad\mbox{for }t\in[t_k,t_{k+1}].
\label{eul}
\end{equation}
It is well known that the order of convergence of the strong error of
discretization is $N^{-1/2}$. Indeed, we have (see~\cite{Ka})
%
%
\begin{equation}
\forall p\geq1, \exists C<+\infty, \forall N\geq1\qquad
\mathbb{E}^{1/p}
\Bigl[\sup_{t\leq T}|X_t-\bar{X}_t|^p
\Bigr]\leq\frac{C}{\sqrt{N}}.
\end{equation}
See Section~\ref{Sec_res_std} for a more precise statement. This
upper-bound gives the correct order of convergence since according to
Remark 3.6 \cite{kurtz5}, when $\sigma$ and $b$ are continuously
differentiable, $(\sqrt{N}(X_t-\bar{X}_t))_{t\leq T}$ converges in law
as $N$ goes to $\infty$ to some diffusion limit which is nonzero as
soon as $\sigma$ is positive and nonconstant (see also \cite
{kurtz6} and \cite{jac1} where stable convergence is also proved). When
$\sigma$ is constant, then the Euler scheme coincides with the Milstein
scheme, and the strong order of convergence is $N^{-1}$.

On the other hand, the order of convergence of the weak error of
discretization is always $N^{-1}$. For example, according to \cite{tt},
when $\sigma$ and $b$ are $C^\infty$ with bounded derivatives of all
orders and $f\dvtx \mathbb{R}\to\mathbb{R}$ is $C^\infty$ with
polynomial growth together with its derivatives, then for each integer
$L\geq1$, the expansion
%
%
\begin{equation}
\mathbb{E} \bigl[f(X_T) \bigr]-\mathbb{E} \bigl[f(
\bar{X}_T) \bigr]=\sum_{l=1}^L
\frac
{a_l}{N^l}+{\mathcal O} \bigl(N^{-(L+1)} \bigr)\label{deverrfaib}
\end{equation}
in powers of $N^{-1}$ holds for the weak error. The bound
$|\mathbb{E}[f(\bar{X}_T)]-\mathbb{E}[f(X_T)]|\leq\frac{C}{N}$ holds
when $\sigma,b$ and $f$ are $C^4$ with the same growth assumptions.
When $f$ is only assumed to be measurable and bounded, it is proved in
\cite{bt1,bt2} that the expansion (\ref{deverrfaib}) still holds for
$L=1$ if $b$ and $\sigma$ are smooth functions satisfying an
hypoellipticity condition. Under uniform ellipticity, Guyon
\cite{gu} even extends this expansion by only assuming that $f$ is a
tempered distribution acting on the densities of both $X_T$ and
$\bar{X}_T$.

In view of financial applications, the weak error analysis gives the
convergence rate to $0$ of the discretization bias introduced when
replacing $X$ by its Euler scheme $\bar{X}$ for the computation of the
price $\mathbb{E}[f(X_T)]$ of a vanilla European option with payoff $f$
and maturity $T$ written on $X$. Let $\mathcal{C}$ denote the space
$C([0,T],\mathbb{R})$ of continuous paths endowed with the sup norm.
When dealing with exotic options with payoff $F\dvtx
\mathcal{C}\to\mathbb{R}$ Lipschitz continuous,
\[
\bigl\llvert\mathbb{E} \bigl[F(X) \bigr]-\mathbb{E} \bigl[F(\bar{X})
\bigr]
\bigr\rrvert\leq\mathbb{E}\bigl|F(X)-F(\bar{X})\bigr|\leq\frac{C}{\sqrt{N}},
\]
where the second inequality follows from the strong error estimate. But
the first inequality is very rough and prevents us from taking
advantage of the cancellations in the mean which occur and permit us to
obtain the upper-bound $\frac{C}{N}$ for vanilla options. The weak
error analysis has been performed for specific path-dependent payoffs,
typically when $F(X)=f(X_T,Y_T)$ with $Y_t$ a function of $(X_s)_{0\leq
s\leq t}$ such that $((X_t,Y_t))_{0\leq t\leq T}$ is a Markov process.
The cases $Y_t=\int_0^tX_s\,ds$ and $Y_t=\max_{0\leq s\leq t}X_s$,
respectively, correspond to Asian \cite{thesestemam} and barrier \cite
{gob1,gob2,gobmen} or lookback options \cite{thess}. But no general
theory has been developped so far to analyze the weak trajectorial
error. The Wasserstein distance between the laws $\mathcal{L}(X)$ and
$\mathcal{L}(\bar{X})$ of $X$ and $\bar{X}$ defined by
\[
\mathcal{W}_1 \bigl(\mathcal{L}(X),\mathcal{L}(\bar{X}) \bigr)=\sup
_{F\dvtx \mathcal{C}\to\mathbb
{R}\dvtx \mathrm{Lip}(F)\leq1}\bigl|\mathbb{E} \bigl[F(\bar{X}) \bigr]-\mathbb{E} \bigl[F(X)
\bigr]\bigr|,
\]
where $\mathrm{Lip}(F)$ denotes the Lipschitz constant of $F$ is the
appropriate measure to deal with the whole class of exotic Lipschitz
payoffs. Notice that this distance has already been used in the context
of discretization schemes for SDEs: in the multidimensional setting, by
a clever rotation of the driving Brownian motion, Cruzeiro, Malliavin
and Thalmaier~\cite{cruz} construct a modified Milstein scheme which
does not involve the simulation of iterated Brownian integrals and with
order of convergence $N^{-1}$ for the Wasserstein distance. A simpler
scheme with the same convergence properties is exhibited in
\cite{js} for usual stochastic volatility models.

The weak and strong error estimations recalled above imply that
%
%
\begin{equation}
\exists c, C<+\infty, \forall N\geq1\qquad \frac{c}{N}\leq\mathcal{W}_1
\bigl(\mathcal{L}(X),\mathcal{L}(\bar{X}) \bigr)\leq\frac{C}{\sqrt{N}}.
\label{minmajw1}
\end{equation}
A very nice feature of the Wasserstein distance is its primal
representation in the Kantorovitch duality theory. This representation
is obtained by choosing $p=1$, $E=\mathcal{C}$ and
$(\mu,\nu)=(\mathcal{L}(X),\mathcal{L}(\bar{X}))$ in the general
definition
%
%
\begin{equation}
\mathcal{W}_p(\mu,\nu)= \biggl(\inf_{\pi\in\Pi(\mu,\nu)}\int
_{E\times E}|x-y|^p\pi(dx,dy) \biggr)^{1/p},
\label{defwas}
\end{equation}
where $p\in[1,+\infty)$, $(E,|~|)$ is a normed vector space, $\mu$ and
$\nu$ are two probability measures on $E$ endowed with its Borel
sigma-field and the infimum is computed on the set $\Pi(\mu,\nu)$ of
probability measures on $E\times E$ with respective marginals $\mu$ and
$\nu$; see, for instance, Remark 6.5 page 95 \cite{villani}. When one
is able to exhibit some coupling $(Y,\bar{Y})$ with
$Y\stackrel{\mathcal{L}}{=}X$ and
$\bar{Y}\stackrel{\mathcal{L}}{=}\bar{X}$, then the law of
$(Y,\bar{Y})$ belongs\vspace*{1pt} to
$\Pi(\mathcal{L}(X),\mathcal{L}(\bar{X}))$ and necessarily
$\mathcal{W}_p(\mathcal{L}(X),\mathcal{L}(\bar{X}))\leq
\mathbb{E}^{1/p} [\sup_{t\in[0,T]}|Y_t-\bar{Y}_t|^p ]$. For the obvious
coupling $(Y,\bar{Y})=(X,\bar{X})$ obtained by choosing the same
driving Brownian motion for the diffusion and its Euler scheme, one
recovers the upper bound in (\ref{minmajw1}) from the strong error
analysis. The main result of the present paper is the construction of a
better coupling which leads to the upper bound
\[
\forall p\geq1, \forall\varepsilon>0, \exists C<+\infty, \forall N\geq1
\qquad
\mathcal{W}_p \bigl(\mathcal{L}(X),\mathcal{L}(\bar{X}) \bigr)\leq
\frac{C}{N^{2/3-\varepsilon}}
\]
proved in Section~\ref{sec_pathwise} under additional regularity
assumptions on the coefficients and uniform ellipticity. To construct
this coupling, we first obtain in Section~\ref{sec_marginal} a
time-uniform estimation of the Wasserstein distance between the
respective laws $\mathcal{L}(X_t)$ and $\mathcal{L}(\bar{X}_t)$ of
$X_t$ and $\bar{X}_t$,
\[
\forall p\geq1, \exists C<+\infty, \forall N\geq1\qquad
\sup_{t\in[0,T]} \mathcal{W}_p \bigl(\mathcal{L}(X_t),\mathcal{L}(
\bar{X}_t) \bigr)\leq\frac{C\sqrt{\log(N)}}{N}.
\]
Previously, in Section~\ref{Sec_res_std}, we recalled well-known
results concerning the moments and the dependence on the initial
condition of the solution to the SDE~(\ref{sde}) and its Euler scheme.
Also, we make explicit the dependence of the strong error estimations
$\mathbb{E} [\sup_{s\le t}|\bar{X}_s-X_s|^p]$ with respect to
$t\in[0,T]$, which will play a key role in our analysis.

\section{Basic estimates on the SDE and its Euler scheme}\label{Sec_res_std}

We recall some well-known results concerning the flow defined
by~(\ref{sde}) (see, e.g., Karatzas and Shreve~\cite{KS}, page 306) and
its Euler approximation.

%
\begin{aprop}
Let us denote by $(X^{x}_t)_{t\in[0,T]}$ the solution of (\ref{sde}),
starting from $x\in\mathbb{R}$. One has that for any $p\geq1$, the existence
of a positive constant $C\equiv C(p,T)$ such that
%
%
\begin{eqnarray}
\forall x \in\mathbb{R}\qquad
\mathbb{E} \Bigl[\sup_{t\in[0,T]}\bigl|X^{x}_t\bigr|^p\Bigr]&\leq& C\bigl(1+|x|\bigr)^p, \label{momenteds}
\\
\qquad\quad \forall x\in\mathbb{R}, \forall s \leq t \leq T
\qquad \mathbb{E} \Bigl[\sup _{u\in[s,t]}\bigl|X^{x}_{u}-X^{x}_{s}\bigr|^p
\Bigr] &\leq& C\bigl(1+|x|\bigr)^p(t-s)^{p/2},\label{accroisseds}
\\
\forall x, y\in\mathbb{R}\qquad \mathbb{E} \Bigl[\sup_{t\in
[0,T]}\bigl|X^{x}_t-X^{y}_t\bigr|^p
\Bigr] &\leq& C|y-x|^p.\label{cieds}
\end{eqnarray}
\end{aprop}

%
\begin{aprop}\label{vitfort_prop}
Let $(\bar{X}^{x}_t)_{t\in[0,T]}$
denote the
Euler scheme~(\ref{eul}) starting from~$x$.
For any $p\in[1,\infty)$, there exists a positive constant $C\equiv
C(p,T)$ such that
%
%
\begin{eqnarray}
\forall N\geq1, \forall x\in\mathbb{R}\qquad \mathbb{E} \Bigl[\sup
_{t\in[0,T]}\bigl|\bar{X}^{x}_t\bigr|^p
\Bigr]&\leq& C\bigl(1+|x|\bigr)^p,\label{momenteul}
\\
\hspace*{30pt}\forall N\geq1, \forall x\in\mathbb{R}, \forall t\in[0,T]\qquad \mathbb{E} \Bigl[\sup_{r\in[0,t]}\bigl|\bar{X}^{x}_r-X^{x}_r\bigr|^p
\Bigr]&\leq&\frac{C t^{p/2}(1+|x|)^p}{N^{p/2}}.\label{vitfort}
\end{eqnarray}
\end{aprop}

The moment bound~(\ref{momenteul}) for the Euler scheme holds in fact
as soon as the drift and the diffusion coefficients have a sublinear
growth. The strong convergence order is established in
Kanagawa~\cite{Ka} for Lipschitz and bounded coefficients. In fact, it
is straightforward to extend Kanagawa's proof to merely Lipschitz
coefficients by using the estimates~(\ref{momenteds})
and~(\ref{momenteul}) and obtain
%
%
\begin{equation}
\label{vitfort2}
\qquad\forall N\geq1, \forall x\in\mathbb{R}, \forall t\in[0,T]\qquad \mathbb{E} \Bigl[\sup_{r\in[0,t]}\bigl|\bar{X}^{x}_r-X^{x}_r\bigr|^p
\Bigr]\leq\frac{C (1+|x|)^p}{N^{p/2}}.
\end{equation}
The estimate~(\ref{vitfort}) precises the dependence on $t$. This
slight improvement will in fact play a crucial role in constructing the
coupling between the diffusion and the Euler scheme. We prove it for
the sake of completeness, even though the arguments are standard.

\begin{pf*}{Proof of~(\ref{vitfort})}
Let $\tau_s=\sup\{t_i, t_i\le s\}$ denote the last discretization time
before~$s$. We have $\bar{X}^{x}_t-X^x_t=\int_0^t
b(\bar{X}^{x}_{\tau_s})-b(X^x_s) \,ds + \int_0^t
\sigma(\bar{X}^{x}_{\tau_s})-\sigma(X^x_s) \,dW_s$. By the Jensen and
Burkholder--Davis--Gundy inequalities,
\begin{eqnarray*}
&& \mathbb{E} \Bigl[ \sup_{r \in[0,t]} \bigl|\bar{X}^{x}_r-X^x_r\bigr|^p
\Bigr]
\\
&&\qquad \le 2^p \biggl( \mathbb{E} \biggl[ \biggl(\int
_0^t \bigl|b \bigl(\bar{X}^{x}_{\tau_s}
\bigr)-b \bigl(X^x_s \bigr)\bigr| \,ds \biggr)^p
\biggr]
\\
&&\hspace*{48pt}{} + C_p\mathbb{E} \biggl[ \biggl(\int_0^t
\bigl(\sigma\bigl(\bar{X}^{x}_{\tau_s} \bigr)-\sigma
\bigl(X^x_s \bigr) \bigr)^2 \,ds
\biggr)^{p/2} \biggr] \biggr)
\\
&&\qquad \le 2^p \biggl( t^{p-1}\int_0^t
\mathbb{E} \bigl[ \bigl|b \bigl(\bar{X}^{x}_{\tau
_s} \bigr)-b
\bigl(X^x_s \bigr)\bigr|^p \bigr] \,ds
\\
&&\hspace*{48pt}{}+ C_pt^{p/2-1}\int_0^t \mathbb{E} \bigl[
\bigl|\sigma\bigl(\bar{X}^{x}_{\tau_s} \bigr)- \sigma\bigl(X^x_s
\bigr)\bigr|^p \bigr]\,ds \biggr).
\end{eqnarray*}

Denoting by $\mathrm{Lip}(\sigma)$ the finite Lipschitz constant of
$\sigma$, we have $|\sigma(\bar{X}^{x}_{\tau_s})-\sigma(X^x_s)|\le
\mathrm{Lip}(\sigma)(|\bar{X}^{x}_{\tau_s}-X^x_{\tau_s}|+|X^x_{\tau_s}-X^x_s|)$.
Thus, (\ref{accroisseds}) and~(\ref{vitfort2}) yield\break
$\mathbb{E}[|\sigma(\bar{X}^{x}_{\tau_s})-\sigma(X^x_s)|^p]\le \frac{C
(1+|x|)^p}{N^{p/2}}, $ and the same bound holds for $b$
replacing~$\sigma$. Since $t^p\leq T^{p/2}t^{p/2}$, we easily conclude.
\end{pf*}

\section{The Wasserstein distance between the marginal laws}\label{sec_marginal}

In this section, we are interested in finding an upper bound for the
Wasserstein distance between the marginal laws of the SDE~(\ref{sde})
and its Euler scheme. It is well known that the optimal coupling
between two one-dimensional random variables is obtained by the inverse
transform sampling. Thus, let $F_t$ and $\bar{F}_t$ denote the
respective cumulative distribution functions of $X_t$ and $\bar{X}_t$.
The $p$-Wasserstein distance between the time-marginals of the
solution\vadjust{\goodbreak}
to the SDE and its Euler scheme is given by (see Theorem~3.1.2
in~\cite{raru})
%
%
\begin{equation}
\mathcal{W}_p \bigl(\mathcal{L}(X_t),\mathcal{L}(
\bar{X}_t) \bigr)= \biggl(\int_0^1\bigl|F_t^{-1}(u)-
\bar{F}_t^{-1}(u)\bigr|^p\,du \biggr)^{1/p}.
\label{wpinv}
\end{equation}
Let us state now the main result of this section. We set
\begin{eqnarray*}
C^k_b&=& \bigl\{f\dvtx \mathbb{R} \rightarrow\mathbb{R}\ k \mbox{
times continuously differentiable  s.t. }
\\
&&\hspace*{115pt}\bigl\|f^{(i)} \bigr\|_\infty<
\infty, 0\le i\le k \bigr\}.
\end{eqnarray*}

\begin{ahyp}\label{hyp_wass_marginal}
Let $a=\sigma^2$. We assume that $a, b \in C^2_b$, $a''$
is globally \mbox{$\gamma$-}H\"older continuous with $\gamma>0$ and
\[
\exists\underline{a}>0, \forall x\in\mathbb{R}, a(x)\geq\underline
{a} \mbox{ (uniform ellipticity)}.
\]
\end{ahyp}

Since $\sigma$ is Lipschitz continuous, under
Hypothesis~\ref{hyp_wass_marginal}, we have either
$\sigma\equiv\sqrt{a}$ or $\sigma\equiv-\sqrt{a}$. From now on, we
assume without loss of generality that $\sigma\equiv\sqrt{a}$ which is
a $C^2_b$ function bounded from below by the positive constant
$\underline{\sigma}=\sqrt{\underline{a}}$.

\begin{theorem}\label{wasun}
Under Hypothesis~\ref{hyp_wass_marginal}, we have for
any $p\ge1$,
\[
\forall N\geq1\qquad \sup_{t\in[0,T]}\mathcal{W}_p \bigl(
\mathcal{L}(X_t),\mathcal{L}(\bar{X}_t) \bigr)\leq
\frac{C\sqrt{\log(N)}}{N},
\]
where $C$ is a positive constant that only depends on $p$, $T$,
$\underline{a}$ and ($\|a^{(i)}\|_\infty$, $\|b^{(i)}\|_\infty$, $0\le
i\le2$) and does not depend on the initial condition~$x\in\mathbb{R}$.
\end{theorem}

\begin{arem}\label{w1unif}
When $p=1$, the slightly better bound
$\sup_{t\in[0,T]}\mathcal{W}_1(\mathcal{L}(X_t),\allowbreak \mathcal{L}(\bar
{X}_t))\leq\frac{C}{N}$ holds if $\sigma$ is uniformly elliptic,
according to~\cite{thesesbai}, Chapter~3. This is proved in a
multidimensional setting for $C^\infty$ coefficients $\sigma$ and $b$
with bounded derivatives by extending the results of \cite{gu} but can
also be derived from a result of Gobet and Labart~\cite{goblab} only
supposing that $b,\sigma\in C^{3}_b$. Let $p_t(x,y)$ and
$\bar{p}_t(x,y)$ denote, respectively, the density of~$X^{0,x}_t$ and
$\bar{X}^{0,x}_t$. Then Theorem~2.3 in~\cite{goblab} gives the
existence of a constant $c>0$ and a finite nondecreasing function $K$
(depending on the upper bounds of $\sigma$ and $b$ and their
derivatives) such that
\[
\forall(t,x,y) \in(0,T]\times\mathbb{R}^2\qquad
\bigl|p_t(x,y)-\bar{p}_t(x,y)\bigr|\leq\frac{TK(T)}{N t}\exp\biggl(-\frac
{c|x-y|^2}{t}\biggr).
\]
As remarked in~\cite{thesesbai}, Chapter~3, for $f\dvtx \mathbb{R}\to
\mathbb{R}$ a Lipschitz continuous function with Lipschitz constant not
greater than one, one deduces that
\begin{eqnarray*}
\bigl|\mathbb{E} \bigl[f(X_t) \bigr]-\mathbb{E} \bigl[f(
\bar{X}_t) \bigr]\bigr|&=& \biggl\llvert\int_{\mathbb
{R}}
\bigl(f(y)-f(x) \bigr) \bigl(p_t(x,y)-\bar{p}_t(x,y)
\bigr)\,dy \biggr\rrvert
\\
&\leq&\frac
{K(T)T}{N t}\int_{\mathbb{R}}|y-x|\exp\biggl(-
\frac{c|x-y|^2}{t} \biggr)\,dy
\\
&=&\frac{K(T)T}{cN},
\end{eqnarray*}
which gives $\sup_{t\leq
T}\mathcal{W}_1(\mathcal{L}(X_t),\mathcal{L}(\bar{X}_t))\leq
\frac{CK(T)T}{N}$ by the dual formulation of the $1$-Wasserstein
distance. 
\end{arem}

Our approach consists of controlling the time evolution of the
Wasserstein distance. To do so, we need to compute the evolution of
both $F_t^{-1}(u)$ and $\bar{F}_t^{-1}(u)$. In the two next
propositions, we derive partial differential equations satisfied by
these functions by integrating in space the Fokker--Planck equations
and then applying the implicit function theorem.

%
\begin{aprop}\label{propevolftm1}
Assume that
Hypothesis~\ref{hyp_wass_marginal} holds.
Then for any $t\in(0,T]$, the cumulative distribution function
$x\mapsto F_t(x)$ is invertible with inverse denoted by $F_t^{-1}(u)$.
Moreover, the function $(t,u)\mapsto F_t^{-1}(u)$ is $C^{1,2}$ on
$(0,T]\times(0,1)$ and satisfies
%
%
\begin{equation}
\partial_t F_t^{-1}(u)=-\frac{1}{2}
\partial_u \biggl(\frac
{a(F_t^{-1}(u))}{\partial_u F_t^{-1}(u)} \biggr)+b \bigl(F_t^{-1}(u)
\bigr).\label{fpinvfr}
\end{equation}
\end{aprop}

%
\begin{aprop}\label{propevolbarftm1}
Assume that $\sigma$ and $b$ have linear growth $\exists C>0$, $\forall
x\in\mathbb{R}$, $|\sigma(x)|+|b(x)|\leq C(1+|x|)$ and that uniform
ellipticity holds, $\exists\underline{a}>0$, $\forall x\in\mathbb{R}$,
$a(x)\geq\underline{a}$. Then for any $t\in(0,T]$, $\bar{X}_t$ admits a
density $\bar{p}_t(x)$ with respect to the Lebesgue measure and its
cumulative distribution function $x\mapsto\bar{F}_t(x)$ is invertible
with inverse denoted by $\bar{F}_t^{-1}(u)$. Moreover, for each
$k\in\{0,\ldots,N-1\}$, the function $(t,u)\mapsto\bar{F}^{-1}_t(u)$ is
$C^{1,2}$ on $(t_k,t_{k+1}]\times(0,1)$ and, on this set, it is a
classical solution of
%
%
\begin{eqnarray}
\label{eqevolbarftm1} 
\partial_t \bar{F}_t^{-1}(u)&=&-
\frac{1}{2}\partial_u \biggl(\frac
{\alpha_t(u)}{\partial_u \bar{F}_t^{-1}(u)} \biggr)+
\beta_t(u),
\end{eqnarray}
where $\alpha_t(u)=\mathbb{E}[a(\bar{X}_{t_k})|\bar{X}_t=\bar
{F}_t^{-1}(u)]$ and
$\beta_t(u)=\mathbb{E}[b(\bar{X}_{t_k})|\bar{X}_t=\bar{F}_t^{-1}(u)]$.
\end{aprop}

The proofs of these two propositions are postponed to
Appendix~\ref{App_Sec1}. Let us mention here that
Proposition~\ref{propevolftm1} also holds when $b'$ is only H\"older
continuous: the Lipschitz assumption on~$b'$ is needed later to prove
Theorem~\ref{wasun}. The PDEs (\ref{fpinvfr})~and~(\ref{eqevolbarftm1})
enable us to compute the time derivative of the $p$th power of the
Wasserstein distance (\ref{wpinv}) and prove, again in
Appendix~\ref{App_Sec1} the following key lemma.

%
\begin{alem}\label{lemmajoderwp}
Under Hypothesis~\ref{hyp_wass_marginal}, for $p\geq2$, the function
$t\mapsto\mathcal{W}_p^p(\mathcal{L}(X_t),\allowbreak
\mathcal{L}(\bar{X}_t))$ is continuous on $[0,T]$, and its first order
distribution\vadjust{\goodbreak} derivative\break
$\partial_t\mathcal{W}_p^p(\mathcal{L}(X_t),\mathcal{L}(\bar{X}_t))$ is
an integrable function on $[0,T]$. Moreover, $dt$ a.e.,
%
%
\begin{eqnarray}\label{majoderwp}
&& \partial_t\mathcal{W}_p^p \bigl(\mathcal{L}(X_t),\mathcal{L}(\bar{X}_t) \bigr)\nonumber
\\[-1pt]
&&\qquad \leq C \biggl(\mathcal{W}_p^p \bigl(\mathcal{L}(X_t),
\mathcal{L}(\bar{X}_t) \bigr)
\nonumber\\[-9pt]\\[-9pt]
&&\hspace*{45pt}{} +\int_0^1\bigl|F_t^{-1}(u)-
\bar{F}_t^{-1}(u)\bigr|^{p-1}\bigl|b \bigl(
\bar{F}_t^{-1}(u) \bigr)-\beta_t(u)\bigr|\,du\nonumber
\\[-1pt]
&&\hspace*{45pt}{}+\int_0^1\bigl|F_t^{-1}(u)-
\bar{F}_t^{-1}(u)\bigr|^{p-2} \bigl(a \bigl(
\bar{F}_t^{-1}(u) \bigr)-\alpha_t(u)
\bigr)^2\,du \biggr),\nonumber
\end{eqnarray}
where $C$ is a positive constant that only depends on $p$,
$\underline{a}$, $\|a'\|_\infty$ and $\|b'\|_\infty$.
\end{alem}

The last ingredient of the proof of Theorem~\ref{wasun} is the next
lemma, the proof of which is also postponed in Appendix~\ref{App_Sec1}.

%
\begin{alem}\label{malcal} Let $\tau_t=\sup\{t_i, t_i\le t\}$ denote
the last discretization time before~$t$. Under
Hypothesis~\ref{hyp_wass_marginal}, we have for all $p\geq1$,
\[
\exists C<+\infty, \forall N\geq1, \forall t\in[0,T]\qquad \mathbb{E} \bigl[
\bigl\llvert\mathbb{E} [ W_{t}-W_{\tau_{t}%
}|\bar{X}_{t} ]
\bigr\rrvert^{p} \bigr] \leq C \biggl(\frac{1}{N\vee(N^2t)}
\biggr)^{p/2}.
\]
\end{alem}

\begin{pf*}{Proof of Theorem~\ref{wasun}}
Since
$\mathcal{W}_p(\mathcal{L}(X_t),\mathcal{L}(\bar{X}_t))\le\mathcal
{W}_{p'}(\mathcal{L}(X_t),\mathcal{L}(\bar{X}_t)) $ for $p\le p'$, it
is enough to prove the estimation for $p\geq2$. Therefore we suppose
without loss of generality that $p\ge2$. Let
$\psi_p(t)=\mathcal{W}^2_p(\mathcal{L}(X_t),\mathcal{L}(\bar{X}_t))$
and\looseness=-1
%
\begin{eqnarray}
\mbox{for any integer } k\geq1\qquad h_k(x)=k^{-2/p}h(kx)\nonumber
\\
\eqntext{\mbox{where }h(x)= \cases{x^{2/p}, &\quad if $x\geq1$,
\vspace*{2pt}\cr
1+\dfrac{2}{p}(x-1), &\quad otherwise.}}
\end{eqnarray}\looseness=0
Since $h_k$ is $C^1$ and nondecreasing, Lemma~\ref{lemmajoderwp} and
H\"older's inequality imply that
\begin{eqnarray*}
&& h_k \bigl(\psi^{p/2}_p(t) \bigr)
\\[-2pt]
&&\qquad = h_k
\bigl(\mathcal{W}_p^p \bigl(\mathcal{L}(X_0),
\mathcal{L}(\bar{X}_0) \bigr) \bigr)+\int_0^th_k'
\bigl(\psi^{p/2}_p(s) \bigr)\partial_s
\mathcal{W}_p^p \bigl(\mathcal{L}(X_s),
\mathcal{L}(\bar{X}_s) \bigr)\,ds
\\[-2pt]
&&\qquad \leq h_k(0)
+C\int_0^th_k'
\bigl(\psi^{p/2}_p(s) \bigr)
\\[-2pt]
&&\hspace*{91pt} {}\times\biggl[\psi ^{p/2}_p(s)
\\[-2pt]
&&\hspace*{107pt}{}
 +\psi^{(p-1)/2}_p(s)
\biggl(\int _0^1\bigl|b \bigl(\bar{F}_s^{-1}(u)
\bigr)-\beta_s(u)\bigr|^p\,du \biggr)^{1/p} \nonumber
\\
&&\hspace*{107pt}
{}+\psi^{(p-2)/2}_p(s) \biggl(\int_0^1\bigl|a
\bigl(\bar{F}_s^{-1}(u) \bigr)-\alpha
_s(u)\bigr|^p\,du \biggr)^{2/p} \biggr]\,ds.
\end{eqnarray*}
Since for fixed $x\geq0$, the sequence $(h'_k(x))_k$ is nondecreasing
and converges to $\frac{2}{p}x^{(2/p)-1}$ as $k\to\infty$, one may take
the limit in this inequality thanks to the monotone convergence theorem
and remark that the image of the Lebesgue measure on $[0,1]$ by
$\bar{F}_s^{-1}$ is the distribution of $\bar{X}_s$ to deduce
%
%
\begin{eqnarray}\label{pregron}
\psi_p(t)&\leq&\frac{2C}{p}\int_0^t
\psi_p(s)+\psi^{1/2}_p(s)\mathbb{E}
^{1/p} \bigl(\bigl|b(\bar{X}_s)-\mathbb{E} \bigl(b(
\bar{X}_{\tau_s})|\bar{X}_s \bigr)\bigr|^p \bigr)
\nonumber\\[-8pt]\\[-8pt]
&&{} + \mathbb{E}^{2/p} \bigl(\bigl|a(\bar{X}_s)-\mathbb{E} \bigl(a(
\bar{X}_{\tau _s})|\bar{X}_s \bigr)\bigr|^p \bigr) \,ds.\nonumber
\end{eqnarray}
One has
\begin{eqnarray*}
a(\bar{X}_{\tau_s})-a(\bar{X}_s)&=&a'(
\bar{X}_s)\sigma(\bar{X}_{s}) (W_{\tau_s}-W_s)
\\[-1pt]
&&{} -a'(\bar{X}_s) \bigl[ \bigl(\sigma(\bar{X}_{\tau
_s})-\sigma(
\bar{X}_{s}) \bigr) (W_s-W_{\tau_s})+b(
\bar{X}_{\tau
_s}) (s-\tau_s) \bigr]
\\[-2pt]
&&{}+(\bar{X}_{\tau_s}-\bar{X}_s)\int_0^1a'
\bigl(v \bar{X}_{\tau_s}+(1-v)\bar{X}_s
\bigr)-a'( \bar{X}_s)\,dv.
\end{eqnarray*}
Using Jensen's inequality, the boundedness assumptions on $a,b$ and
their derivatives and Lemma \ref{malcal}, one gets
\begin{eqnarray*}
&& \mathbb{E} \bigl(\bigl|a(\bar{X}_s)-
\mathbb{E} \bigl(a( \bar{X}_{\tau_s})|\bar{X}_s
\bigr)\bigr|^p \bigr)
\\
&&\qquad \leq C\mathbb{E} \bigl(\bigl|\sigma a'(
\bar{X}_s)\bigr|^p\bigl| \mathbb{E}\bigl((W_s-W_{\tau_s})|
\bar{X}_s \bigr)\bigr|^p \bigr)
\\
&&\quad\qquad{}+C\mathbb{E} \bigl((s-\tau_s)^p+\bigl| \bigl(\sigma(\bar
{X}_{\tau
_s})-\sigma(\bar{X}_{s}) \bigr) (W_s-W_{\tau_s})\bigr|^p+|
\bar{X}_{\tau
_s}-\bar{X}_s|^{2p} \bigr)
\\
&&\qquad \leq\frac{C}{N^{p/2}\vee(N^ps^{p/2})}.
\end{eqnarray*}
The same bound holds with $a$ replaced by $b$. With (\ref{pregron}) and
Young's inequality, one deduces
\begin{eqnarray*}
\psi_p(t)&\leq& C\int_0^t
\psi_p(s)+\frac{\psi^{1/2}_p(s)}{\sqrt{N}\vee(N\sqrt{s})}+\frac
{1}{N\vee(N^2s)}\,ds
\\
&\leq& C\int _0^t\psi_p(s)+
\frac{1}{N\vee(N^2s)}\,ds.
\end{eqnarray*}
One concludes by Gronwall's lemma.
\end{pf*}

%
\begin{arem} When $a(x)\equiv a$ is constant, the term $\mathbb
{E}^{2/p}
(|a(\bar{X}_s)-\break\mathbb{E}(a(\bar{X}_{\tau_s})|\bar{X}_s)|^p )$ in
(\ref{pregron}) vanishes and the above reasoning ensures that $\bar
{\psi}_p(t)$ defined as $\sup_{s\in[0,T]}\psi_p(s)$ satisfies
\begin{eqnarray*}
\bar{\psi}_p(t)&\leq& C\int_0^t \bar{\psi}_p(s) \,ds
+C\bar{\psi}^{1/2}_p(t) \int
_0^t\frac{1}{\sqrt{N}\vee(N\sqrt{s})}\,ds
\\
&\le& C\int _0^t\bar{\psi}_p(s)
\,ds + \frac{1}{2}\bar{\psi}_p(t) +\frac{C^2 (T+1)^2}{2N}.
\end{eqnarray*}
By Gronwall's lemma, we recover the estimation
$\sup_{t\in[0,T]}\mathcal{W}_p(\mathcal{L}(X_t),\break \mathcal{L}(\bar
{X}_t))\leq\frac{C}{N}$
which is also a consequence of the strong order of convergence of the
Euler scheme when the diffusion coefficient is constant.
\end{arem}

\section{The Wasserstein distance between the pathwise laws}\label{sec_pathwise}
We now state the main result of the paper.
%
%
\begin{ahyp}\label{hyp_wass_pathwise}
We assume that $a \in C^4_b, b \in C^3_b$ and
\[
\exists\underline{a}>0, \forall x\in\mathbb{R}\qquad a(x)\geq\underline{a}
\mbox{ (uniform ellipticity)}.
\]
%
\end{ahyp}
Clearly, Hypothesis~\ref{hyp_wass_pathwise} implies
Hypothesis~\ref{hyp_wass_marginal}.

\begin{theorem}\label{main_thm}
Under Hypothesis~\ref{hyp_wass_pathwise}, we have
\[
\forall p\geq1, \forall\varepsilon>0, \exists C<+\infty, \forall N\geq1
\qquad \mathcal{W}_p \bigl(\mathcal{L}(X),\mathcal{L}(\bar{X}) \bigr)\leq
\frac{C}{N^{2/3-\varepsilon}}.
\]
\end{theorem}

Before proving the theorem, let us state some of its consequences for
the pricing of lookback options. It is well known (see, e.g.,
\cite{glas} page 367) that if $(U_k)_{0\leq k\leq N-1}$ are independent
random variables uniformly distributed on $[0,1]$ and independent from
the Brownian increments $(W_{t_{k+1}}-W_{t_{k}})_{0\leq k\leq N-1}$
then
$\bar{\hspace*{-1.2pt}\bar{X}}\stackrel{\mathrm{def}}{=}\frac{1}{2}\max_{0\leq
k\leq N-1} (\bar{X}_{t_k}+\bar{X}_{t_{k+1}}+\sqrt{(\bar
{X}_{t_{k+1}}\,{-}\,\bar{X}_{t_k})^2-2\sigma^2(\bar{X}_{t_k})t_1\ln(U_k)} ) $
is such\vspace*{-2pt} that $
(\bar{X}_0,\bar{X}_{t_1},\ldots,\bar{X}_{T},\bar{\hspace*{-1.2pt}\bar{X}}
)\stackrel{\mathcal{L}}{=}(\bar{X}_0,\bar{X}_{t_1},\ldots,\bar
{X}_{T},\max_{t\in [0,T]}\bar{X}_t)$.
%
%
\begin{acor}
If $f\dvtx \mathbb{R}^2\to\mathbb{R}$ is Lipschitz continuous, then,
under Hypothesis~\ref{hyp_wass_pathwise},
%
%
\begin{eqnarray}\label{vitfaiblook}
&& \forall\varepsilon>0, \exists C<+\infty, \forall N\geq1
\nonumber\\[-6pt]\\[-12pt]
&&\qquad \Bigl\llvert \mathbb{E} \Bigl[f \Bigl(X_T,\max_{t\in[0,T]}X_t
\Bigr) \Bigr]-\mathbb{E}
\bigl[f(\bar{X}_T,\bar{\hspace*{-1.2pt}\bar{X}}) \bigr]
\Bigr\rrvert\leq\frac{C}{N^{2/3-\varepsilon}}.\nonumber
\end{eqnarray}
\end{acor}

To our knowledge, this result appears to be new. Of course, when $f$ is
also differentiable with respect to its second variable, one has
\begin{eqnarray*}
&& \mathbb{E} \Bigl[f \Bigl(X_T,\max_{t\in[0,T]}X_t \Bigr)
\Bigr]
\\
&&\qquad =\mathbb{E} \bigl[f (X_T,x_0 ) \bigr]+
\int_{x_0}^{+\infty}\mathbb{E} \bigl[
\partial_2f(X_T,x)1_{\{\max_{t\in[0,T]}X_t\geq x\}} \bigr]\,dx.
\end{eqnarray*}
One could try to combine the weak error analysis for the first term on
the right-hand side with Theorem 2.3 \cite{gob2} devoted to barrier
options to obtain the order $N^{-1}$ instead on $N^{-2/3+\varepsilon}$
in (\ref{vitfaiblook}). Unfortunately, one cannot succeed for two main
reasons. First, it is not clear whether the estimation in Theorem 2.3
\cite{gob2} is preserved by integration over $[x_0,+\infty)$. More
importantly, for this estimation to hold, a structure condition on the
payoff function implying that $\partial_2f(x,x)=0$ for all $x\geq x_0$
is needed.

\begin{pf*}{Proof of Theorem \ref{main_thm}}
We first deduce from Theorem \ref{wasun} some bound on the Wasserstein
distance between the finite-dimensional marginals of the diffusion $X$
and its Euler scheme $\bar{X}$ on a coarse time-grid. For
$m\in\{1,\ldots,N-1\}$, we set $n=\lfloor N/m\rfloor$ and define
\[
s_l=\frac{lmT}{N}\qquad\mbox{for } l\in\{0,\ldots,n-1\}\mbox{ and } s_n=T.
\]
We will use this coarse time-grid $(s_l)_{1\leq l\leq n}$ to
approximate the supremum norm on $\mathcal{C}$ and therefore we endow
consistently $\mathbb{R}^n$ with the norm $|(x_1,\ldots,\break x_n)|=\max
_{1\leq l\leq n}|x_l|$. Combining the next proposition, the proof of
which is postponed in Appendix~\ref{App_Sec2} with Theorem \ref{wasun},
one obtains that
%
%
\begin{equation}
\mathcal{W}_p \bigl(\mathcal{L}(X_{s_1},\ldots,X_{s_n}),\mathcal{L}(\bar{X}_{s_1},\ldots,
\bar{X}_{s_n}) \bigr)\leq\frac{C\sqrt{\log
N}}{m},\label{wasfindim}
\end{equation}
where the constant $C$ does not depend on $(m,N)$.
%
%
\begin{aprop}\label{prop_wass_multi}
Let $\mathbb{R}^n$ be endowed
with the norm $|(x_1,\ldots,x_n)|=\max_{1\leq l\leq n}|x_l|$. For any
$p\geq1$, there is a constant $C$ not depending on $n$ such that
\[
\mathcal{W}_p \bigl(\mathcal{L}(X_{s_1},\ldots,X_{s_n}),\mathcal{L}(\bar{X}_{s_1},\ldots,
\bar{X}_{s_n}) \bigr)\leq Cn\sup_{0\leq t\leq T,x\in
\mathbb{R}}
\mathcal{W}_p \bigl(\mathcal{L} \bigl(\bar{X}^{x}_t
\bigr),\mathcal{L} \bigl(X^{x}_t \bigr) \bigr).
\]
\end{aprop}
There is a probability measure
$\pi(dx_1,\ldots,dx_n,d\bar{x}_1,\ldots,d\bar{x}_n)$ in $\Pi
(\mathcal{L}(X_{s_1},\allowbreak\ldots,
X_{s_n}),\mathcal{L}(\bar{X}_{s_1},\ldots,\bar{X}_{s_n}))$ which
attains the Wasserstein distance in the left-hand side of
(\ref{wasfindim}); see, for instance, Theorem 3.3.11 \cite{raru},
according to which $\pi$ is the law of
$(X_{s_1},\ldots,X_{s_n},\xi_{s_1},\ldots,\xi_{s_n})$ with
$(\xi_{s_1},\ldots,\xi_{s_n})\in\partial_{|~|}\varphi
(X_{s_1},\ldots,\break X_{s_n})$ where $\partial_{|~|}\varphi$ is the
subdifferential, for the above defined norm $|~|$ on $\mathbb{R}^n$, of
some \mbox{$|~|$-}convex function $\varphi$. Let
$\tilde{\pi}(x_1,\ldots,x_n,d\bar{x}_1,\ldots,d\bar{x}_n)$ denote a
regular conditional probability of $(\bar{x}_1,\ldots,\bar{x}_n)$ given
$(x_1,\ldots,x_n)$ when $\mathbb{R}^{2n}$ is endowed with $\pi$ and
$(\bar{Y}_{s_1},\ldots,\bar{Y}_{s_n})$ be distributed according to
$\tilde{\pi}(X_{s_1},\ldots,X_{s_n},\break d\bar{x}_1,\ldots,d\bar{x}_n)$. The
vector $(X_{s_1},\ldots,X_{s_n},\bar{Y}_{s_1},\ldots,\bar {Y}_{s_n})$
is distributed according to $\pi$ so that
%
%
\begin{eqnarray}\label{xybar}
(\bar{Y}_{s_1},\ldots,\bar{Y}_{s_n})&\stackrel{
\mathcal{L}} {=}&(\bar{X}_{s_1},\ldots,\bar{X}_{s_n})\quad\mbox{and}
\nonumber\\[-8pt]\\[-8pt]
\mathbb{E}^{1/p} \Bigl[\max_{1\leq l\leq n}|X_{s_l}-
\bar{Y}_{s_l}|^p \Bigr]&\leq&\frac{C\sqrt{\log N}}{m}.\nonumber
\end{eqnarray}
Let $p_t(x,y)$ denote the transition density of the SDE~(\ref{sde}) and
$\ell_t(x,y)=\log(p_t(x,y))$. According to Appendix \ref{diff_bridge}
devoted to diffusion bridges, the processes
\[
\biggl(W^l_t=\int
_{s_l}^t \bigl(dW_s-
\sigma(X_s)\partial_x\ell_{s_{l+1}-s}(X_s,X_{s_{l+1}})
\,ds \bigr),t \in[s_l,s_{l+1}) \biggr) _{0\leq l\leq n-1}
\]
are independent Brownian motions independent from
$(X_{s_1},\ldots,X_{s_n})$. We suppose from now on that the vector
$(\bar{Y}_{s_1},\ldots,\bar{Y}_{s_n})$ has been generated independently
from these processes and so will be all the random variables and
processes needed in the remaining of the proof (see in particular the
construction of $\beta$ below). Moreover, again by Appendix
\ref{diff_bridge}, the solution of
%
%
\begin{equation}\label{defz}
\cases{
\displaystyle Z^{x,y}_t=x+
\int_{s_l}^t\sigma\bigl(Z^{x,y}_s
\bigr)\,dW^l_s
\vspace*{2pt}\cr
\displaystyle\hspace*{31pt}{} +\int_{s_l}^t \bigl[b \bigl(Z^{x,y}_s \bigr)+\sigma^2 \bigl(Z^{x,y}_s
\bigr)\partial_x\ell _{s_{l+1}-s} \bigl(Z^{x,y}_s,y \bigr) \bigr]\,ds,
\vspace*{2pt}\cr
\hspace*{200pt} t\in[s_l,s_{l+1}),
\vspace*{4pt}\cr Z^{x,y}_{s_{l+1}}=y}
\end{equation}
is distributed according to the conditional law of
$(X_t)_{t\in[s_l,s_{l+1}]}$ given $(X_{s_l},\break X_{s_{l+1}})=(x,y)$ and for
each $l\in\{0,\ldots,n-1\}$, one has $(Z^{X_{s_l},X_{s_{l+1}}}_t)_{t\in
[s_l,s_{l+1}]}=(X_t)_{t\in[s_l,s_{l+1}]}$.

In order to construct a good coupling between $\mathcal{L}(X)$ and
$\mathcal{L}(\bar{X})$, a natural idea would be to extend
$(\bar{Y}_{s_1},\ldots,\bar{Y}_{s_n})$ to a process
$(\bar{Y}_t)_{t\in[0,T]}$ with law $\mathcal{L}(\bar{X})$ by defining
for each $l\in\{0,\ldots,n-1\}$, $(\bar{Y}_t)_{t\in[s_l,s_{l+1}]}$ as
the process obtained by inserting the Brownian motion $W^l$, the
starting point $\bar{Y}_{s_l}$ and the ending point $\bar{Y}_{s_{l+1}}$
in the It\^{o}'s decomposition of the conditional dynamics of
$(\bar{X}_t)_{t\in[s_l,s_{l+1}]}$ given $\bar{X}_{s_l}=x$ and
$\bar{X}_{s_{l+1}}=y$. Unfortunately, even if this Euler scheme bridge
is deduced by a simple transformation of the Brownian bridge on a
single time-step, it becomes a complicated process\vspace*{1pt} when
the difference between the starting and ending times is larger than
$\frac{T}{N}$ because of the lack of Markov property. At the end of the
proof, we will choose the difference $s_{l+1}-s_l$ of order
$\frac{T}{N^{1/3}}$ and, therefore,\vspace*{-2pt} much larger than the
time-step $\frac{T}{N}$. In addition, it is not clear how to compare
the paths of the diffusion bridge and the Euler scheme bridge driven by
the same Brownian motion. That is, why we are going to introduce some
new process $(\tilde{\chi}_t)_{t\in[0,T]}$ such that the comparison
will be performed at the diffusion bridge level, which is not so
easy~yet.

To construct $\tilde{\chi}$, we are going to exhibit a Brownian motion
$(\beta_t)_{t\in[0,T]}$ such that $\bar{Y}_{s_1},\ldots,\bar{Y}_{s_n}$
are the values on the coarse time-grid $(s_l)_{1\leq l\leq n}$ of the
Euler scheme~(\ref{eul}) driven by $\beta$ instead of $W$. The
extension $(\bar{Y}_t)_{t\in[0,T]}$ with law $\mathcal{L}(\bar{X})$ is
then simply defined as the whole Euler scheme driven by $\beta$:
%
\begin{eqnarray}
\bar{Y}_t=\bar{Y}_{t_k}+\sigma(\bar{Y}_{t_k}) (
\beta_t-\beta_{t_k})+b(\bar{Y}_{t_k})
(t-t_k),\nonumber
\\
\eqntext{ t\in[t_{k},t_{k+1}], 0\le k\le N-1.}
\end{eqnarray}
The construction of $\beta$ is postponed at the end of the
present proof. One then defines
\[
\chi_t=\bar{Y}_{s_l}+\int_{s_{l}}^t
\sigma(\chi_{s})\,d\beta_s+\int_{s_{l}}^tb(
\chi_s)\,ds,\qquad t\in[s_{l},s_{l+1}), 0\le l\le
n-1.
\]

Notice that the process $\chi=(\chi_t)_{t\in[0,T]}$ which evolves
according to the SDE~(\ref{sde}) with $\beta$ replacing $W$ on each
time-interval $[s_l,s_{l+1})$ is c\`adl\`ag: discontinuities may arise
at the points $\{s_{l+1}, 0\leq l\leq n-1\}$. We denote by
$\chi_{s_{l+1}-}$ its left-hand limit at time $s_{l+1}$ and set
$\chi_{T}=\chi_{s_n-}$. The strong error estimation (\ref{vitfort})
will permit us to estimate the difference between the processes
$\bar{Y}$ and $\chi$. For the subsequent choice of $\beta$, we do not
expect the processes $\chi$ and $X$ to be close. Nevertheless, the
process $\tilde{\chi}$ obtained by setting
\[
\forall l\in\{0,\ldots,n-1\}, \forall t\in[s_l,s_{l+1})
\qquad
\tilde{\chi}_t=Z^{\chi_{s_l},\chi
_{s_{l+1}-}}_t\quad\mbox{and}\quad\tilde{\chi}_T=\chi_T,
\]
where $Z^{x,y}$ is defined in (\ref{defz}) is such that
$\mathcal{L}(\tilde{\chi})=\mathcal{L}(\chi)$ by
Propositions~\ref{prop_bridge} and~\ref{prop_bridge2}. On each coarse
time-interval $[s_l,s_{l+1})$ the diffusion bridges associated with $X$
and $\tilde{\chi}$ are driven by the same Brownian motion $W^l$.
Moreover the differences $|X_{s_l}-Y_{s_l}|$ between the starting
points and $|X_{s_{l+1}}-\chi_{s_{l+1}-}|\leq
|X_{s_{l+1}}-Y_{s_{l+1}}|+| Y_{s_{l+1}}-\chi_{s_{l+1}-}|$ between the
ending points is controlled by (\ref{xybar}) and the above mentionned
strong error estimation. That is, why one may expect to obtain a good
estimation of the difference between the processes $X$ and
$\tilde{\chi}$. By the triangle inequality and since
$\mathcal{L}(\bar{X})=\mathcal{L}(\bar{Y})$ and
$\mathcal{L}(\tilde{\chi})=\mathcal{L}(\chi)$,
%
%
\begin{eqnarray}\label{triangle}
\qquad\mathcal{W}_p \bigl(\mathcal{L}(\bar{X}),\mathcal{L}(X) \bigr)&\leq&
\mathcal{W}_p \bigl(\mathcal{L}(\bar{X}),\mathcal{L}(\chi) \bigr)+
\mathcal{W}_p \bigl(\mathcal{L}(\chi),\mathcal{L}(X) \bigr)
\nonumber
\\
&\leq&\mathbb{E}^{1/p} \Bigl[\sup_{t\in[0,T]}|
\bar{Y}_t-\chi_t|^p \Bigr]+\mathbb{E}
^{1/p} \Bigl[\sup_{t\in[0,T]}|X_t-\tilde{
\chi}_t|^p \Bigr],
\end{eqnarray}
where, for the definition of
$\mathcal{W}_p(\mathcal{L}(\bar{X}),\mathcal{L}(\chi))$ and
$\mathcal{W}_p(\mathcal{L}(\chi),\mathcal{L}(X))$, the space of
c\`adl\`ag sample-paths from $[0,T]$ to $\mathbb{R}$ is endowed with
the supremum norm. Let us first estimate the first term on the
right-hand side. From~(\ref{vitfort}), we get
\[
\mathbb{E} \Bigl[\sup_{t\in[s_l,s_{l+1})}|\bar{Y}_t-\chi
_t|^{p} \big|\bar{Y}_{s_l} \Bigr]\leq C
\frac{m^{p/2}(1+|\bar{Y}_{s_l}|)^{p}}{N^{p}},
\]
where the constant $C$ does not depend on $(N,m)$. We deduce that
\begin{eqnarray*}
\mathbb{E} \Bigl[\sup_{t\in[0,T]}|\bar{Y}_t-
\chi_t|^{p} \Bigr]
&=&\mathbb{E} \Bigl[\max _{0\leq l\leq n-1}\sup_{t\in[s_l,s_{l+1})}|\bar{Y}_t-\chi
_t|^{p} \Bigr]
\\
&\leq& \sum_{l=0}^{n-1}
\mathbb{E} \Bigl[\mathbb{E} \Bigl[\sup_{t\in
[s_l,s_{l+1})}|
\bar{Y}_t-\chi_t|^{p} \big|\bar{Y}_{s_l}
\Bigr] \Bigr]
\\
&\leq& C\frac{m^{p/2}}{N^{p}}\sum_{l=0}^{n-1}
\mathbb{E} \bigl[\bigl(1+|\bar{Y}_{s_l}|\bigr)^{p} \bigr]
\\
&\leq& C
\frac{m^{p/2-1}}{N^{p-1}},
\end{eqnarray*}
where we used (\ref{momenteul}) for the last inequality. As a
consequence,
%
%
\begin{equation}
\label{Wass_chi_se} \mathbb{E}^{1/p} \Bigl[
\sup_{t\in
[0,T]}|\bar{Y}_t-\chi_t|^p
\Bigr] \le C\frac{m^{1/2-1/p}}{N^{1-1/p}}.
\end{equation}

Let us now estimate the second term on the right-hand side of
(\ref{triangle}). By Proposition~\ref{prop_bridge2} and since for
$l\in\{0,\ldots,n-1\}$, $\chi_{s_l}=\bar{Y}_{s_l}$,
\begin{eqnarray*}
\sup_{t\leq T}|X_t-\tilde{\chi}_t|&=&\max _{0\leq l\leq
n-1}\sup_{t\in
[s_l,s_{l+1})}\bigl|Z^{X_{s_l},X_{s_{l+1}}}_t-Z^{\chi_{s_l},\chi
_{s_{l+1}-}}_t\bigr|
\\
&\leq& C\max_{0\leq l\leq n-1}|X_{s_l}-\bar{Y}_{s_l}| \vee|X_{s_{l+1}}-\chi_{s_{l+1}-}|.
\end{eqnarray*}
Since, by the triangle inequality and the continuity of $\bar{Y}$,
\begin{eqnarray*}
|X_{s_{l+1}}-\chi_{s_{l+1}-}|&\leq&|X_{s_{l+1}}-\bar
{Y}_{s_{l+1}}|+|\bar{Y}_{s_{l+1}}-\chi_{s_{l+1}-}|
\\
&\leq&
|X_{s_{l+1}}-\bar{Y}_{s_{l+1}}|+\sup_{t\in[0,T]}|
\bar{Y}_{t}-\chi_{t}|,
\end{eqnarray*}
one deduces that
\[
\sup_{t\leq T}|X_t-\tilde{\chi}_t|\leq C
\Bigl(\max_{1\leq l\leq
n}|X_{s_l}-\bar{Y}_{s_l}|+
\sup_{t\in[0,T]}|\bar{Y}_{t}-\chi_{t}| \Bigr).
\]
Combined with (\ref{xybar}) and (\ref{Wass_chi_se}), this implies
\begin{eqnarray*}
\mathbb{E}^{1/p} \Bigl[\sup_{t\leq T}|X_t- \tilde{\chi}_t|^p
\Bigr]&\leq& C\mathbb{E} ^{1/p} \Bigl[\max_{1\leq l\leq
n}|X_{s_l}-\bar{Y}_{s_l}|^p \Bigr]+C\mathbb{E}^{1/p}
\Bigl[\sup_{t\in[0,T]}| \bar{Y}_{t}-\chi_{t}|^p \Bigr]
\\
&\leq& C \biggl(
\frac{\sqrt{\log N}}{m}+\frac{m^{1/2-1/p}}{N^{1-1/p}} \biggr).
\end{eqnarray*}
Plugging this inequality together with (\ref{Wass_chi_se}) in
(\ref{triangle}), we deduce that
\[
\mathcal{W}_p \bigl(\mathcal{L}(X),\mathcal{L}(\bar{X}) \bigr) \leq
C \biggl(\frac{\sqrt{\log N}}{m}+\frac{m^{1/2-1/p}}{N^{1-1/p}} \biggr)
\]
and conclude by choosing $m=\lfloor N^{2/3} \rfloor$ that for
$p\geq\frac{1}{3\varepsilon}$,
$\mathcal{W}_p(\mathcal{L}(X),\mathcal{L}(\bar{X}))\leq\frac{C}{N^{2/3-\varepsilon}}$.
When $\frac{1}{3\varepsilon}>1$, the conclusion follows for
$p\in[1,\frac{1}{3\varepsilon})$ since
$\mathcal{W}_p(\mathcal{L}(X),\allowbreak\mathcal{L}(\bar{X}))\leq\mathcal
{W}_{1/3\varepsilon}(\mathcal{L}(X),\mathcal{L}(\bar{X}))$.

To complete the proof, we still have to construct the Brownian motion
$\beta$. We first reconstruct on the fine time grid $(t_k)_{1\leq k\leq
N}$ an Euler scheme $(\bar{Y}_{t_k},0\le k\le N)$ interpolating the
values on the coarse grid $(s_l)_{1\leq l\leq n}$. Let us denote by
$\bar{p}(x,y)$ the density of the law
$\mathcal{N}(x+b(x)T/N,\sigma(x)^2T/N)$ of the Euler scheme starting
from~$x$ after one time step~$T/N$. Thanks to the ellipticity
assumption, we have $\bar{p}(x,y)>0$ for any $x,y\in\mathbb{R}$.
Conditionally on $(\bar{Y}_{s_1},\ldots,\bar{Y}_{s_n})$, we generate
independent random vectors
\[
(\bar{Y}_{s_{l-1}+t_1},\ldots,\bar{Y}_{s_{l-1}+t_{m-1}})_{1\leq
l\leq n-1}\quad\mbox{and}\quad (\bar{Y}_{s_{n-1}+t_1},\ldots,\bar{Y}_{t_{N-1}})
\]
with respective densities
\[
\frac{\bar{p}(\bar{Y}_{s_{l-1}},x_1)\bar{p}(x_1,x_2)\cdots\bar
{p}(x_{n-1},\bar{Y}_{s_{l}})
}{\int_{\mathbb{R}^{n-1}}\bar{p}(\bar{Y}_{s_{l-1}},y_1)\bar
{p}(y_1,y_2)\cdots\bar{p}(y_{n-1},\bar{Y}_{s_{l}})\,dy_1\cdots
dy_{n-1}}
\]
and
\[
\frac{\bar{p}(\bar{Y}_{s_{n-1}},x_1)\bar{p}(x_1,x_2)\cdots \bar
{p}(x_{N-1-m(n-1)},\bar{Y}_{s_{n}}) }{\int_{\mathbb
{R}^{N-1-m(n-1)}}\bar{p}(\bar{Y}_{s_{n-1}},y_1)\bar{p}(y_1,y_2)\cdots
\bar{p}(y_{N-1-m(n-1)},\bar{Y}_{s_{n}})\,dy_1\cdots dy_{N-1-m(n-1)}}.
\]
This ensures that $(\bar{Y}_{t_k})_{0\leq k\leq
n}\stackrel{\mathcal{L}}{=}(\bar{X}_{t_k})_{0\leq k\leq n}$. Then we
get, thanks to the ellipticity condition, that
$ (\frac{1}{\sigma(\bar{Y}_{t_{k-1}})}(\bar{Y}_{t_k}-\bar
{Y}_{t_{k-1}}-b(\bar{Y}_{t_{k-1}})) )_{1\le k\le N}$~are
independent centered Gaussian variables with variance~$T/N$. By using
independent Brownian bridges, we can then construct a Brownian motion
$(\beta_t)_{t\in[0,T]}$ such that
\[
\beta_{t_k}-\beta_{t_{k-1}}= \frac{1}{\sigma(\bar{Y}_{t_{k-1}})} \bigl(
\bar{Y}_{t_k}-\bar{Y}_{t_{k-1}}-b(\bar{Y}_{t_{k-1}})
\bigr),
\]
which completes the construction.
\end{pf*}

\section*{Conclusion}
In this paper, we prove that the order of convergence of the
Wasserstein distance $\mathcal{W}_p$ on the space of continuous paths
between the laws of a uniformly elliptic one-dimensional diffusion and
its Euler scheme with $N$-steps is not worse that
$N^{-2/3+\varepsilon}$. In view of a possible extension to
multidimensional settings, two main difficulties have to be resolved.
First, we take advantage of the optimality of the inverse transform
coupling in dimension one to obtain a uniform bound on the Wasserstein
distance between the marginal laws with optimal rate $N^{-1}$ up to a
logarithmic factor. In dimension $d>1$, the optimal coupling between
two probability measures on $\mathbb{R}^d$ is not available, which
makes the estimation of the Wasserstein distance between the marginal
laws much more complicated even if, for $\mathcal{W}_1$, the order
$N^{-1}$ may be deduced from the results of \cite{goblab}; see Remark
\ref{w1unif}. Next, one has to generalize the estimation on diffusion
bridges given by Proposition~\ref{prop_bridge2} which we deduce from
the Lamperti transform in dimension $d=1$.

In the perspective of the multi-level Monte Carlo method introduced by
Giles~\cite{giles}, coupling with order of convergence
$N^{-2/3+\varepsilon}$ the Euler schemes with $N$ and $2N$ steps would
also be of great interest for variance reduction, especially in
multidimensional situations where the Milstein scheme is not feasible;
see \cite{js} for the implementation of this idea in the example of a
discretization scheme devoted to usual stochastic volatility models.
But this does not seem obvious from our nonconstructive coupling
between the Euler scheme and its diffusion limit. For both the
derivation of the order of convergence of the Wasserstein distance on
the path space and the explicitation of the coupling, the limiting step
in our approach is Proposition \ref{prop_wass_multi}. In this
proposition, we bound the dual formulation of the Wasserstein distance
between $n$-dimensional marginals by the Wasserstein distance between
one-dimensional marginals multiplied by $n$.

Even if the order of convergence of the Wasserstein distance on the
path space obtained in the present paper may not be optimal, it
provides the first significant step from the order $N^{1/2}$ obtained
with the trivial coupling where the diffusion and the Euler scheme are
driven by the same Brownian motion.
%
%
\begin{appendix}
\section{Proofs of Section~\lowercase{\protect\texorpdfstring{\ref{sec_marginal}}{2}}}\label{App_Sec1}

\begin{pf*}{Proof of Proposition \ref{propevolftm1}}
According to \cite{friedman}, Theorems 5.4 and 4.7, for any
$t\in(0,T]$, the solution $X_t$ of (\ref{sde}) starting from $X_0=x_0$
admits a density $p_t(x)$ w.r.t. the Lebesgue measure on the real line,
the function $(t,x)\mapsto p_t(x)$ is $C^{1,2}$ on
$(0,T]\times\mathbb{R}$, and on this set,\vspace*{-1pt} it is a
classical solution of the Fokker--Planck equation
%
%
\begin{equation}
\partial_t p_t(x)=\tfrac{1}{2}\partial
_{xx} \bigl(a(x)p_t(x) \bigr)-\partial_x
\bigl(b(x)p_t(x) \bigr). \label{fp}
\end{equation}
Moreover, the following Gaussian bounds hold:
%
%
\begin{eqnarray}\label{gb}
\bigl|p_t(x)\bigr|+\sqrt{t}\bigl|\partial_x p_t(x)\bigr|\leq\frac{C}{\sqrt{t}}e^{-(x-x_0)^2/Ct}
\nonumber\\[-10pt]\\[-10pt]
\eqntext{\exists C>0, \forall t\in(0,T], \forall x\in\mathbb{R}.}
\end{eqnarray}
The partial derivatives $\partial_x F_t(x)=p_t(x)$ and
$\partial_{xx}F_t(x)=\partial_xp_t(x)$ exist and are continuous on
$(0,T]\times\mathbb{R}$. For $0<s<t\leq T$ and $y\leq x$, integrating
(\ref{fp}) over $[s,t]\times[y,x]$, then letting $y\to-\infty$ thanks
to (\ref{gb}), one obtains
$F_t(x)-F_s(x)=\int_{s}^t\frac{1}{2}\partial_{x}(a(x)p_r(x))-b(x)p_r(x)\,dr$.
By continuity of the integrand w.r.t. $(r,x)$ one deduces that the
partial derivative $\partial_tF_t(x)$ exists and is continuous on
$(0,T]\times\mathbb{R}$. So, $(t,x)\mapsto F_t(x)$ is $C^{1,2}$ on
$(0,T]\times\mathbb{R}$ and solves
%
%
\begin{equation}
\label{fpfr}\partial_t F_t(x)=\tfrac{1}{2}
\partial_{x} \bigl(a(x)\partial_x F_t(x)
\bigr)-b(x)\partial_x F_t(x).
\end{equation}

According to Aronson \cite{aron}, the density is also bounded from
below by some Gaussian kernel $\exists c>0, \forall
(t,x)\in(0,T]\times\mathbb{R}, |p_t(x)|\geq
\frac{c}{\sqrt{t}}e^{-(x-x_0)^2/ct}$. This enables us to apply the
implicit function theorem to $(t,x,u)\mapsto F_t(x)-u$ to deduce that
the inverse $u\mapsto F_t^{-1}(u)$ of $x\mapsto F_t(x)$ is $C^{1,2}$ in
the variables $(t,u)\in(0,T]\times(0,1)$ and solves
\begin{eqnarray*}
\partial_t F_t^{-1}(u)&=&-\frac{\partial_t F_t}{\partial_x
F_t}
\bigl(F_t^{-1}(u) \bigr)
\\
&=&-\frac{1}{2}\partial_{x} \bigl(a(x)\partial_x
F_t(x) \bigr)\big|_{x=F_t^{-1}(u)}\partial_u
F_t^{-1}(u)+b \bigl(F_t^{-1}(u)
\bigr)
\\
&=&-\frac{1}{2}\partial_{u} \biggl(\frac{a(F_t^{-1}(u))}{\partial_u
F_t^{-1}(u)}
\biggr)+b \bigl(F_t^{-1}(u) \bigr),
\end{eqnarray*}
where we used\vspace*{-1pt} (\ref{fpfr}) for the second equality and
$\partial_u F_t^{-1}(u)=\frac{1}{\partial_xF_t( F_t^{-1}(u))}$ for both
the second and the third equalities.
\end{pf*}

\begin{pf*}{Proof of Proposition \ref{propevolbarftm1}}
For $t\in(0,t_1]$,
$\bar{X_t}$ admits the Gaussian density with mean $x_0+b(x_0)t$ and
variance $a(x_0)t$. By induction on $k$ and independence of
$W_t-W_{t_k}$ and $\bar{X}_{t_k}$ in (\ref{eul}), one checks that for
$k\in\{1,\ldots,\allowbreak n-1\}$, $\bar{X}_{t_k}$ admits a density
$\bar{p}_{t_k}(x)$ and that for $t\in(t_{k},t_{k+1}]$, $(\bar
{X}_{t_k},\bar{X_t})$ admits the density
\[
\rho(t_k,t,y,x)=\bar{p}_{t_k}(y)\frac{\exp({-
{(x-y-b(y)(t-t_k))^2}/{2a(y)(t-t_k)}})}{\sqrt{2\pi a(y)(t-t_k)}}.
\]
The marginal density $\bar{p}_t(x)=\int_\mathbb{R}\bar
{p}_{t_k}(y)\frac{\exp({-({x-y-b(y)(t-t_k)^2})/{2a(y)(t-t_k)}})}{\sqrt{2\pi
a(y)(t-t_k)}}\,dy$ of $\bar{X}_t$ is continuous on
$(t_k,t_{k+1}]\times\mathbb{R}$ by Lebesgue's theorem and positive.

Let $N(x)=\int_{-\infty}^xe^{-y^2/2}\frac{dy}{\sqrt{2\pi}}$ denote the
cumulative distribution function of the standard Gaussian law and
$k\in\{0,\ldots,N-1\}$. Again by the independence structure in
(\ref{eul}), for $(t,x)\in(t_k,t_{k+1}]\times\mathbb{R}$,
$\bar{F}_t(x)=\break \mathbb{E} (N (\frac{x-\bar{X}_{t_k}-b(\bar
{X}_{t_k})(t-t_k)}{\sqrt{a(\bar{X}_{t_k})(t-t_k)}} ) )$. One~has
\begin{eqnarray*}
\partial_t N \biggl(\frac{x-y-b(y)(t-t_k)}{\sqrt{a(y)(t-t_k)}} \biggr
)&=& - \biggl(\frac{x-y-b(y)(t-t_k)}{2\sqrt{2\pi a(y)(t-t_k)^3}}+\frac
{b(y)}{\sqrt{2\pi a(y)(t-t_k)}} \biggr)
\\
&&{}\times \exp\biggl(-\frac
{(x-y-b(y)(t-t_k))^2}{2a(y)(t-t_k)}\biggr).
\end{eqnarray*}
By the growth assumption on $\sigma$ and $b$, one easily checks that
$\forall k\in\{0,\ldots,N\}$, $\mathbb{E}(\bar{X}^2_{t_k})<+\infty$.
With the uniform ellipticity assumption, one deduces by a standard
uniform integrability argument that $\bar{F}_t(x)$ is differentiable
w.r.t. $t$ with partial\vadjust{\goodbreak} derivative
%
%
\begin{eqnarray}\label{evolfbart}
\qquad \partial_t\bar{F}_t(x)&=&-\mathbb{E} \biggl[ \biggl(
\frac{x-\bar{X}_{t_k}-b(\bar{X}_{t_k})(t-t_k)}{2\sqrt{2\pi a(\bar
{X}_{t_k})(t-t_k)^3}}+\frac{b(\bar{X}_{t_k})}{\sqrt{2\pi a(\bar
{X}_{t_k})(t-t_k)}} \biggr)
\nonumber\\[-8pt]\\[-8pt]
&&\hspace*{70pt}
{}\times
\exp \biggl(-\frac{(x-\bar{X}_{t_k}-b(\bar
{X}_{t_k})(t-t_k))^2}{2a(\bar{X}_{t_k})(t-t_k)}\biggr) \biggr]\nonumber
\end{eqnarray}
continuous in $(t,x)\in(t_k,t_{k+1}]\times\mathbb{R}$. In the same way,
one checks smoothness of $\bar{F}_t(x)$ in the spatial variable $x$ and
obtains that this function is $C^{1,2}$ on \mbox{$(t_k,t_{k+1}]\times
\mathbb{R}$.}

When $k\geq1$,
\begin{eqnarray*}
&& \mathbb{E} \biggl[b(\bar{X}_{t_k})\frac{\exp (-{(x-\bar
{X}_{t_k}-b(\bar{X}_{t_k})(t-t_k))^2}/{(2a(\bar{X}_{t_k})(t-t_k)}))}{\sqrt{2\pi a(\bar
{X}_{t_k})(t-t_k)}} \biggr]
\\
&&\qquad =\int_{\mathbb{R}}b(y)\rho(t_k,t,y,x)\,dy
\\
&&\qquad =\mathbb{E} \bigl[b(
\bar{X}_{t_k})|\bar{X}_t=x \bigr]\bar{p}_t(x).
\end{eqnarray*}
For $k=0$, even if $(\bar{X}_0,\bar{X}_t)$ has no density, the
equality between the opposite sides of this equation remains true.

Combining Lebesgue's theorem and a similar reasoning, one checks that
\begin{eqnarray*}
&& -\mathbb{E} \biggl[\frac{x-\bar{X}_{t_k}-b(\bar
{X}_{t_k})(t-t_k)}{\sqrt{2\pi a(\bar{X}_{t_k})(t-t_k)^3}}\exp\biggl({-\frac
{(x-\bar{X}_{t_k}-b(\bar
{X}_{t_k})(t-t_k))^2}{2a(\bar{X}_{t_k})(t-t_k)}}\biggr) \biggr]
\\
&&\qquad =\partial
_x\mathbb{E} \biggl[a(\bar{X}_{t_k})\frac{\exp({-{(x-\bar
{X}_{t_k}-b(\bar
{X}_{t_k})(t-t_k))^2}/({2a(\bar{X}_{t_k})(t-t_k)})})}{\sqrt{2\pi a(\bar
{X}_{t_k})(t-t_k)}}
\biggr]
\\
&&\qquad =\partial_x \bigl[\mathbb{E} \bigl(a(\bar{X}_{t_k})|
\bar{X}_t=x \bigr)\bar{p}_t(x) \bigr].
\end{eqnarray*}
With (\ref{evolfbart}), one deduces that
%
%
\begin{equation}
\label{gyongyfbart}
\qquad\partial_t\bar{F}_t(x)=
\tfrac{1}{2}\partial_x \bigl(\mathbb{E} \bigl[a(\bar
{X}_{t_k})|\bar{X}_t=x \bigr]\partial_x
\bar{F}_t(x) \bigr)-\mathbb{E} \bigl[b(\bar{X}_{t_k})|
\bar{X}_t=x \bigr]\partial_x\bar{F}_t(x).
\end{equation}
One checks that the function $(t,u)\mapsto\bar{F}_t^{-1}(u)$ is smooth
and satisfies the partial differential equation (\ref{eqevolbarftm1})
by arguments similar to the ones given at the end of the proof of
Proposition \ref{propevolftm1}.
\end{pf*}

%
\begin{arem}
In the same way, for $k\in\{0,\ldots,N-1\}$, one could prove that on
$(t_k,t_{k+1}]\times\mathbb{R}$, $(t,x)\mapsto\bar{p}_t(x)$ is
$C^{1,2}$ and satisfies the partial differential
\[
\partial_t\bar{p}_t(x)=\tfrac{1}{2}
\partial_{xx} \bigl(\mathbb{E} \bigl[a(\bar{X}_{t_k})|
\bar{X}_t=x \bigr]\bar{p}_t(x) \bigr)-
\partial_x \bigl(\mathbb{E} \bigl[b(\bar{X}_{t_k})|
\bar{X}_t=x \bigr] \bar{p}_t(x) \bigr)\vadjust{\goodbreak}
\]
obtained by spatial derivation of (\ref{gyongyfbart}). This shows that
$(\bar{X}_t)_{t\in[0,T]}$ has the same marginal distributions as the
diffusion process with coefficients given by the above conditional
expectations, which is also a consequence of \cite{gyongy}.
\end{arem}

\begin{pf*}{Proof of Lemma \ref{lemmajoderwp}}
By the continuity of
the paths of $X$ and $\bar{X}$ and the finiteness of
$\mathbb{E} [\sup_{t\leq T}(|X_t|^{p+1}+|\bar
{X}_t|^{p+1}) ]$, one easily checks that
$t\mapsto\mathcal{W}_p^p(\mathcal{L}(X_t),\mathcal{L}(\bar{X}_t))$ is
continuous.

Let $k\in\{0,\ldots,N-1\}$ and $s,t\in(t_k,t_{k+1}]$ with $s\leq t$.
Combining Propositions \ref{propevolftm1}~and~\ref{propevolbarftm1}
with a spatial integration by parts, one obtains for $\varepsilon\in
(0,1/2)$,
%
%
\begin{eqnarray}\label{preipp}
\hspace*{12pt}&& \int_\varepsilon^{1-\varepsilon}\bigl|F_t^{-1}(u)-
\bar{F}_t^{-1}(u)\bigr|^p\,du\nonumber
\\[-2pt]
&&\qquad =\int_\varepsilon^{1-\varepsilon}\bigl|F_s^{-1}(u)-\bar
{F}_s^{-1}(u)\bigr|^p\,du\nonumber
\\[-2pt]
&&\quad\qquad{} +p\int_s^t\int_\varepsilon
^{1-\varepsilon}\bigl|F_r^{-1}(u)-\bar{F}_r^{-1}(u)\bigr|^{p-2}
\bigl(F_r^{-1}(u)-\bar{F}_r^{-1}(u)\bigr)\nonumber
\\[-2pt]
&&\hspace*{92pt}
{}\times\bigl(b \bigl(F_r^{-1}(u) \bigr)-
\beta_r(u) \bigr)\,du\,dr\nonumber
\\[-2pt]
&&\quad\qquad{} +\frac{p(p-1)}{2}\int_s^t\int
_\varepsilon^{1-\varepsilon
}\bigl|F_r^{-1}(u)-
\bar{F}_r^{-1}(u)\bigr|^{p-2} \bigl(
\partial_u F_r^{-1}(u)-\partial_u
\bar{F}_r^{-1}(u) \bigr)
\nonumber\\[-9pt]\\[-9pt]
&&\hspace*{128pt}
{}\times \biggl(\frac
{a(F_r^{-1}(u))}{\partial_u F_r^{-1}(u)}-
\frac{\alpha_r(u)}{\partial
_u \bar{F}_r^{-1}(u)} \biggr)\,du\,dr\nonumber
\\[-2pt]
&&\quad\qquad{}+\frac{p}{2}\int_s^t\bigl|F_r^{-1}(1-
\varepsilon)-\bar{F}_r^{-1}(1-\varepsilon)\bigr|^{p-2}
\bigl(F_r^{-1}(1-\varepsilon)-\bar{F}_r^{-1}(1-
\varepsilon) \bigr)\nonumber
\\[-2pt]
&&\hspace*{68pt}
{}\times\biggl(\frac{\alpha_r(1-\varepsilon
)}{\partial_u \bar{F}_r^{-1}(1-\varepsilon)}-\frac
{a(F_r^{-1}(1-\varepsilon))}{\partial_u F_r^{-1}(1-\varepsilon
)} \biggr)\,dr\nonumber
\\[-2pt]
&&\quad\qquad{} -\frac{p}{2}\int_s^t\bigl|F_r^{-1}(
\varepsilon)-\bar{F}_r^{-1}(\varepsilon)\bigr|^{p-2}
\bigl(F_r^{-1}(\varepsilon)-\bar{F}_r^{-1}(
\varepsilon) \bigr)\nonumber
\\[-2pt]
&&\hspace*{68pt}
{}\times \biggl(\frac{\alpha_r(\varepsilon
)}{\partial_u \bar{F}_r^{-1}(\varepsilon)}-\frac
{a(F_r^{-1}(\varepsilon))}{\partial_u F_r^{-1}(\varepsilon)} \biggr)\,dr.\nonumber
\end{eqnarray}

We are now going to take the limit as $\varepsilon\to0$. We will check
at the end of the proof that
%
%
\begin{eqnarray}\label{IPP_terms}
&& \lim_{u\rightarrow0^+\ \mathrm{or}\ 1^-}\ \sup
_{r\in[s,t]}\frac{a(F_t^{-1}(u))}{\partial_u
F_t^{-1}(u) } \bigl|F_t^{-1}(u)-
\bar{F}_t^{-1}(u) \bigr|^{p-1}
\nonumber\\[-8pt]\\[-8pt]
&&\quad{} + \sup _{r\in [s,t]}\frac{\alpha_t(u) }{\partial_u \bar{F}_t^{-1}(u) }
\bigl|F_t^{-1}(u)- \bar{F}_t^{-1}(u) \bigr|^{p-1}=0,\nonumber
\end{eqnarray}
which enables us to get rid of the two last boundary terms.\vadjust{\goodbreak}

Combining Young's inequality with the uniform ellipticity assumption
and the positivity of $\partial_uF_t^{-1}(u)$ and $\partial_u
\bar{F}_t^{-1}(u)$, one obtains
\begin{eqnarray*}
&& \bigl(\partial_u F_r^{-1}(u)-
\partial_u \bar{F}_r^{-1}(u) \bigr) \biggl(
\frac{a(F_r^{-1}(u))}{\partial_u F_r^{-1}(u)}-\frac{\alpha_r(u)}{\partial
_u \bar{F}_r^{-1}(u)} \biggr)
\\
&&\qquad = \bigl(a \bigl(F_r^{-1}(u) \bigr)-
\alpha_r(u) \bigr)\frac{\partial_u
F_r^{-1}(u)-\partial_u \bar{F}_r^{-1}(u)}{\partial_u F_r^{-1}(u)\vee
\partial_u \bar{F}_r^{-1}(u)}
\\
&&\quad\qquad{} -a \bigl(F_r^{-1}(u) \bigr)\frac{((\partial_u
\bar{F}_r^{-1}(u)-\partial_u F_r^{-1}(u))^+)^2}{\partial_u
F_r^{-1}(u)\partial_u \bar{F}_r^{-1}(u)}
\\
&&\quad\qquad{}-\alpha_r(u)\frac
{((\partial_u F_r^{-1}(u)-\partial_u \bar
{F}_r^{-1}(u))^+)^2}{\partial_u F_r^{-1}(u)\partial_u \bar
{F}_r^{-1}(u)}
\\
&&\qquad\leq\frac{1}{4\underline{a}} \bigl(a \bigl(F_r^{-1}(u)
\bigr)- \alpha_r(u) \bigr)^2+\underline{a}
\frac{(\partial_u F_r^{-1}(u)-\partial
_u \bar{F}_r^{-1}(u))^2}{(\partial_u F_r^{-1}(u)\vee\partial_u \bar
{F}_r^{-1}(u))^2}
\\
&&\quad\qquad{} - \bigl(a \bigl(F_r^{-1}(u) \bigr)\wedge
\alpha_r(u) \bigr)\frac{(\partial
_u \bar{F}_r^{-1}(u)-\partial_u F_r^{-1}(u))^2}{\partial_u
F_r^{-1}(u)\partial_u \bar{F}_r^{-1}(u)}
\\
&&\qquad\leq\frac{1}{4\underline{a}} \bigl(a \bigl(F_r^{-1}(u)
\bigr)- \alpha_r(u) \bigr)^2.
\end{eqnarray*}
Hence, up to the factor $\frac{p(p-1)}{2}$, the third term on the
right-hand side of (\ref{preipp}) is equal to
\begin{eqnarray*}
&&\int_s^t\int_\varepsilon^{1-\varepsilon}\bigl|F_r^{-1}(u)-
\bar{F}_r^{-1}(u)\bigr|^{p-2}
\\[-2pt]
&&\hspace*{38pt}{}\times \biggl[ \bigl(
\partial_u F_r^{-1}(u)-\partial_u
\bar{F}_r^{-1}(u) \bigr) \biggl(\frac{a(F_r^{-1}(u))}{\partial_u
F_r^{-1}(u)}-
\frac{\alpha_r(u)}{\partial_u \bar{F}_r^{-1}(u)} \biggr)
\\
&&\hspace*{174pt}{} -\frac{ (a(F_r^{-1}(u))-\alpha_r(u) )^2}{4\underline
{a}} \biggr]\,du\,dr
\\[-2pt]
&&\quad{} +\frac{1}{4\underline{a}}\int_s^t
\int_\varepsilon^{1-\varepsilon}\bigl|F_r^{-1}(u)-
\bar{F}_r^{-1}(u)\bigr|^{p-2} \bigl(a
\bigl(F_r^{-1}(u) \bigr)-\alpha_r(u)
\bigr)^2\,du\,dr,
\end{eqnarray*}
where the integrand in the first integral is nonpositive. Since
\begin{eqnarray*}
\hspace*{-5pt}&&\int_s^t\hspace*{-1pt}\int_0^1\bigl|F_r^{-1}(u)-
\bar{F}_r^{-1}(u)\bigr|^{p-2}
\\[-1pt]
\hspace*{-5pt}&&\hspace*{25pt}{} \times\hspace*{-0.3pt}  \bigl(\bigl|F_r^{-1\hspace*{-0.3pt}}(u)-\hspace*{-0.2pt}
\bar{F}_r^{-1\hspace*{-0.3pt}}(u)\bigr|\bigl|b \bigl(F_r^{-1\hspace*{-0.3pt}}(u)
\bigr)-\beta_r(u)\bigr|
+ \bigl(a \bigl(F_r^{-1\hspace*{-0.3pt}}(u)
\bigr)-\alpha_r(u) \bigr)^2 \bigr)\,du\,dr
\\[-1pt]
\hspace*{-5pt}&&\qquad \leq2\|b\|_\infty\int_s^t
\mathcal{W}_{p}^{p-1} \bigl(\mathcal{L}(X_r),
\mathcal{L}(\bar{X}_r) \bigr)\,dr
\\[-1pt]
\hspace*{-5pt}&&\quad\qquad{}+4\|a\|^2_\infty
\int_s^t\mathcal{W}_{p}^{p-2}
\bigl(\mathcal{L}(X_r),\mathcal{L}(\bar{X}_r)
\bigr)\,dr<+\infty,
\end{eqnarray*}
one can take the
limit $\varepsilon\to0$ in (\ref{preipp}) using Lebesgue's theorem for
the second term on the right-hand side and combining Lebesgue's theorem
with monotone convergence for the third term to obtain
%
%
\begin{eqnarray}\label{vipp}
&& \mathcal{W}_{p}^{p} \bigl(\mathcal{L}(X_t),
\mathcal{L}(\bar{X}_t) \bigr)\nonumber
\\
&&\qquad =\mathcal{W}_{p}^{p}
\bigl(\mathcal{L}(X_s),\mathcal{L}(\bar{X}_s) \bigr)
\nonumber
\\
&&\quad\qquad{}+p\int_s^t\int_0^{1}\bigl|F_r^{-1}(u)-
\bar{F}_r^{-1}(u)\bigr|^{p-2} \bigl(F_r^{-1}(u)-
\bar{F}_r^{-1}(u) \bigr)
\nonumber\\[-8pt]\\[-8pt]
&&\hspace*{82pt}{}\times \bigl(b \bigl(F_r^{-1}(u)
\bigr)-\beta_r(u) \bigr)\,du\,dr
\nonumber
\\
&&\quad\qquad{}+\frac{p(p-1)}{2}\int_s^t\int
_0^{1}\bigl|F_r^{-1}(u)-\bar
{F}_r^{-1}(u)\bigr|^{p-2} \bigl(\partial_u
F_r^{-1}(u)-\partial_u \bar
{F}_r^{-1}(u) \bigr)\nonumber
\\
&&\hspace*{118pt}
{}\times \biggl(\frac{a(F_r^{-1}(u))}{\partial_u
F_r^{-1}(u)}-
\frac{\alpha_r(u)}{\partial_u \bar{F}_r^{-1}(u)} \biggr)\,du\,dr.\nonumber
\end{eqnarray}
The last term which belongs to $[-\infty,+\infty)$ is finite since so
are all the other terms. We deduce the integrability of
\begin{eqnarray*}
(r,u)&\mapsto&\bigl|F_r^{-1}(u)-\bar{F}_r^{-1}(u)\bigr|^{p-2}
\bigl(\partial_u F_r^{-1}(u)-
\partial_u \bar{F}_r^{-1}(u) \bigr)
\\
&&\times{} \biggl(\frac{a(F_r^{-1}(u))}{\partial_u
F_r^{-1}(u)}-\frac{\alpha_r(u)}{\partial_u \bar{F}_r^{-1}(u)} \biggr)
\end{eqnarray*}
on $[s,t]\times(0,1)$. Similar arguments show that the integrability
property and (\ref{vipp}) remain true for $s=t_k$. By summation, they
remain true for $0\leq s\leq t\leq T$. So the integrability holds on
$[0,T]$ for the derivative in the distributional sense
\begin{eqnarray*}
&& \partial_t\mathcal{W}_{p}^{p} \bigl(
\mathcal{L}(X_t),\mathcal{L}(\bar{X}_t) \bigr)
\\
&&\qquad =p\int_0^{1}\bigl|F_t^{-1}(u)-
\bar{F}_t^{-1}(u)\bigr|^{p-2} \bigl(F_t^{-1}(u)-
\bar{F}_t^{-1}(u) \bigr) \bigl(b \bigl(F_t^{-1}(u)
\bigr)-\beta_t(u) \bigr)\,du
\\
&&\quad\qquad{}+\frac{p(p-1)}{2}\int_0^{1}\bigl|F_t^{-1}(u)-
\bar{F}_t^{-1}(u)\bigr|^{p-2} \bigl(
\partial_u F_t^{-1}(u)-\partial_u
\bar{F}_t^{-1}(u) \bigr)
\\
&&\hspace*{103pt}
{}\times \biggl(\frac
{a(F_t^{-1}(u))}{\partial_u F_t^{-1}(u)}-
\frac{\alpha_t(u)}{\partial
_u \bar{F}_t^{-1}(u)} \biggr)\,du
\\
&&\qquad \leq p\int_0^{1}\bigl|F_t^{-1}(u)-
\bar{F}_t^{-1}(u)\bigr|^{p-2}
\\
&&\hspace*{53pt}
{}\times \biggl[
\bigl(F_t^{-1}(u)-\bar{F}_t^{-1}(u)
\bigr) \bigl(b \bigl(F_t^{-1}(u) \bigr)-
\beta_t(u) \bigr)
\\
&&\hspace*{106pt}{} +\frac{(p-1) (a(F_t^{-1}(u))-\alpha_t(u) )^2}{8\underline {a}} \biggr]\,du.
\end{eqnarray*}
Equation (\ref{majoderwp}) follows by remarking that
\begin{eqnarray*}
&& \bigl(a \bigl(F_t^{-1}(u) \bigr)-\alpha_t(u)
\bigr)^2
\\
&&\qquad \leq2 \bigl(\bigl\|a'\bigr\|_\infty
^2\bigl|F_t^{-1}(u)-\bar{F}_t^{-1}(u)\bigr|^2+
\bigl(a \bigl(\bar{F}_t^{-1}(u) \bigr)-
\alpha_t(u) \bigr)^2 \bigr)
\end{eqnarray*}
and using a similar idea for $|b(F_t^{-1}(u))-\beta_t(u)|$.

To prove (\ref{IPP_terms}) for $0<s\leq t\leq T$, we use the Aronson
estimates recalled in the proof of Proposition \ref{propevolftm1} for
$X_t$ and deduced from Theorem 2.1 \cite{lemairemenozzi}, for the Euler
scheme
%
%
\begin{eqnarray}\label{aronson}
&& \frac{c}{\sqrt{r}}\exp\biggl(-\frac{(x-x_0)^2}{cr} \biggr)
\nonumber\\[-8pt]\\[-8pt]
&&\qquad \le
p_r(x)\wedge\bar{p}_r(x)\le p_r(x)\vee
\bar{p}_r(x)
\le \frac{C}{\sqrt{r}}\exp\biggl(-\frac{(x-x_0)^2}{C r}\biggr).\nonumber
\end{eqnarray}

Setting $K_1=\frac{c}{\sqrt{t}}$, $c_1=cs/2$, $K_2=\frac{C}{\sqrt{s}}$
and $c_2=Ct/2$, one has
%
\begin{eqnarray}\label{aronson2}
K_1 \exp\biggl(-\frac{(x-x_0)^2}{2c_1 } \biggr) \le\rho_r(x) \le
K_2 \exp\biggl(-\frac{(x-x_0)^2}{2c_2 } \biggr)
\nonumber\\[-8pt]\\[-8pt]
\eqntext{\forall r\in[s,t], \forall x \in\mathbb{R},}
\end{eqnarray}
where $\rho_r$ denotes either $p_r$ or $\bar{p}_r$. The four limits in
(\ref{IPP_terms}) can be obtained similarly, and we focus on the one of
$\sup_{r\in[s,t]}\frac{a(F_r^{-1}(u))}{\partial_u
F_r^{-1}(u) } |F_r^{-1}(u)-\bar{F}_r^{-1}(u) |^{p-1}$.
Up to modifying $K_1>0$ and decreasing $c_1>0$, we get
from~(\ref{aronson2}) that
%
\begin{eqnarray}
K_1(x_0-x)\exp\biggl(-\frac{(x-x_0)^2}{2c_1 } \biggr) \le\rho_r(x) \le
K_2 (x_0-x) \exp\biggl(-\frac
{(x-x_0)^2}{2c_2 } \biggr)\nonumber
\\
\eqntext{\forall r\in[s,t], \forall x \le x_0-1,}
\end{eqnarray}
which leads to
\[
\forall x \le x_0-1\qquad K_1c_1\exp\biggl(-
\frac{(x-x_0)^2}{2c_1 } \biggr) \le G_r(x) \le K_2c_2
\exp\biggl(-\frac{(x-x_0)^2}{2c_2 } \biggr),
\]
where $G_r$ denotes either $F_r$ or $\bar{F}_r$. Thus, the inverse
function satisfies
%
%
\begin{equation}
x_0- \sqrt{-2c_2 \log\biggl(\frac{u}{K_2c_2}
\biggr)} \le\bar{F}_r^{-1}(u)\le x_0-
\sqrt{-2c_1 \log\biggl(\frac{u}{K_1c_1} \biggr)}\label{continvrep}
\end{equation}
for $u$ small
enough. The two last inequalities imply that when $x\rightarrow-
\infty$,
\[
\forall r\in[s,t]\qquad \bar{F}_r^{-1} \bigl(F_r(x)
\bigr) \ge x_0 - \sqrt{-2 c_2 \biggl[ \log\biggl(
\frac{K_1c_1}{K_2c_2} \biggr) -\frac{(x-x_0)^2}{2c_1} \biggr]}
\]
and $\sup_{r\in[s,t]}|x-\bar{F}_r^{-1}(F_r(x))
|\mathop{=}\limits_{x\rightarrow- \infty}O(x)$. With the boundedness
of $a$\vspace*{-1pt} and (\ref{aronson2}), we easily deduce that
\[
\lim_{x\rightarrow- \infty} \sup_{r\in[s,t]}a(x)
p_r(x) \bigl|x-\bar{F}_r^{-1}
\bigl(F_r(x) \bigr) \bigr|^{p-1}=0.
\]
Since, by (\ref{continvrep}), $\bar{F}_r^{-1}(u)$
converges to $-\infty$ uniformly in $r\in[s,t]$ as $u$ tends to $0$,
we conclude that
\[
\lim_{u\rightarrow0^+}\sup_{r\in[s,t]}\frac
{a(F_r^{-1}(u))}{\partial_u
F_r^{-1}(u) }
\bigl|F_r^{-1}(u)-\bar{F}_r^{-1}(u)
\bigr|^{p-1}=0.
\]\upqed
\end{pf*}

\begin{pf*}{Proof of Lemma \ref{malcal}}
By Jensen's inequality,
\begin{eqnarray*}
\mathbb{E} \bigl[\bigl|\mathbb{E}(W_t-W_{\tau_t}|
\bar{X}_t)\bigr|^p \bigr]&\leq&\mathbb{E} \bigl[|W_t-W_{\tau_t}|^{p}
\bigr]\leq\frac{C}{N^{p/2}}.
\end{eqnarray*}
Let us now check that the left-hand side is also smaller than
$\frac{C}{t^{p/2}N^{p}}$. To do this, we will study
\[
{\mathbb{E}} \bigl[ (W_{t}-W_{\tau_{t}})g(\bar{X}_{t})
\bigr],
\]
where $g$ is any smooth real valued function.

In order to continue, we need to do various estimations on the Euler
scheme and its Malliavin derivative, which we denote by $D_u\bar{X}_t$.
Let $\eta_{t}=\min\{t_{i};t\leq t_{i}\}$ denote the discretization
time just after $t$. We have $D_u\bar{X}_t=0$ for $u>t$, and
%
\begin{eqnarray}
D_{u}\bar{X}_{t}&=&1_{\{t\leq
\eta_u\}}\sigma(\bar{X}_{\tau_{t}})\nonumber
\\
&&{} +1_{\{t>\eta_u\}} \bigl(1+\sigma'(
\bar{X}_{\tau_{t}}) (W_t-W_{\tau_t})+b'(\bar{X}_{\tau_{t}}) (t-\tau_t) \bigr)D_u
\bar{X}_{\tau_{t}}\nonumber
\\
\eqntext{\mbox{for }u\leq t.}
\end{eqnarray}

Then by induction, one clearly obtains that for $u\le t$,
\begin{eqnarray*}
D_{u}\bar{X}_{t} & =&\sigma(\bar{X}_{\tau_{u}})
\bar{\mathcal{E}}_{u,t},
\\
\bar{\mathcal{E}}_{u,t} & =& \cases{ 1, &\quad if $\tau_{t}
\leq\eta_{u}$,
\vspace*{3pt}
\cr
\bigl( 1+b^{\prime}(
\bar{X}_{\tau_{t}}) (t-\tau_{t})+\sigma^{\prime}(
\bar{X}_{\tau_{t}}) (W_{t}-W_{\tau_{t}}) \bigr), &\quad if $
\eta_{u}=\tau_{t}$,
\vspace*{4pt}
\cr
\displaystyle\prod
_{i=N\eta_{u}/T}^{N\tau_{t}/T-1} \bigl( 1+b^{\prime}(
\bar{X}_{t_{i}}) (t_{i+1}-t_{i})+\sigma^{\prime}(\bar{X}_{t_{i}}) (W_{t_{i+1}}-W_{t_{i}})
\bigr)
\cr
\hspace*{33pt}{}\times\bigl( 1+b^{\prime}(\bar{X}_{\tau_{t}}) (t-\tau_{t})+
\sigma^{\prime}(\bar{X}_{\tau_{t}}) (W_{t}-W_{\tau_{t}})
\bigr),&\quad if $\eta_{u}<\tau_{t}$.}
\end{eqnarray*}

Note that $\bar{\mathcal{E}}$ satisfies the following properties: (1)
$\bar{\mathcal{E}}_{u,t}=\bar{\mathcal{E}}_{\eta(u),t}$ and\break (2)~$\bar{\mathcal{E}}_{t_i,t_j}\bar{\mathcal{E}}_{t_j,t}=\bar{\mathcal
{E}}_{t_i,t}$ for $t_i\le t_j\le t$. We also introduce the process
$\mathcal{E}$ defined by
\[
\mathcal{E}_{u,t}=\exp\biggl( \int_{u}^{t}b^{\prime}(X_{s})-
\frac
{1}{2}%
\sigma^{\prime}(X_{s})^{2}\,ds+
\int_{u}^{t}\sigma^{\prime
}(X_{s})\,dW_{s}
\biggr).
\]
The next lemma, the proof of which is postponed at the end of the
present proof states some useful properties of the processes
$\mathcal{E}$ and $\bar{\mathcal{E}}$.

\begin{alem}\label{lemme_majorations} Let us assume that $b,\sigma\in
C^2_b$. Then we have
%
%
\begin{eqnarray}
\sup_{0\leq s\leq t \le T}{\mathbb{E}} \bigl[
\mathcal{E}_{s,t}^{-p} \bigr]+{\mathbb{E}} \bigl[
\mathcal{E}_{s,t}^{p} \bigr] &\leq& C, \label{eq:propE}
\\
\sup_{0\leq s\leq t \le T}{\mathbb{E}} \bigl[ \bar{\mathcal{E}}_{s,t}^{p}
\bigr] &\leq& C,\label{eqA.13}
\\
\sup_{0\leq s,u\leq t \le T}{\mathbb{E}} \bigl[ |
D_u\bar{\mathcal{E}}_{s,t}|^p+|
D_u \mathcal{E}_{s,t}|^p \bigr] &\leq& C,
\label{eq:A13}
\\
\sup_{0\leq t \le T}{\mathbb{E}} \bigl[ \llvert\mathcal
{E}_{0,t}-\bar{\mathcal{E}}_{0,t}\rrvert
^{p} \bigr] &\leq&\frac
{C}{N^{p/2}}, \label{vitesse_forte}
\end{eqnarray}
where $C$ is a positive constant depending only on $p$ and $T$.
\end{alem}

We next define the localization given by
\[
\psi=\varphi\bigl( \mathcal{E}_{0,t}^{-1} (
\mathcal{E}_{0,t}%
-\bar{\mathcal{E}}_{0,t} ) \bigr).
\]
Here $\varphi\dvtx \mathbb{R\rightarrow}[0,1]$ is a~$C^\infty$
symmetric function so that
\[
\varphi(x)=\cases{ 0, &\quad if $|x|>\frac{1}{2}$, \vspace*{2pt}
\cr
1, &
\quad if $|x|<\frac{1}{4}$.}
\]
One has
\begin{eqnarray*}
\mathbb{E} \bigl[ (W_{t}-W_{\tau_{t}})g(\bar{X}_{t})
\bigr] &=&\mathbb{E} \bigl[ (W_{t}-W_{\tau_{t}})g(
\bar{X}_{t})\psi\bigr]+\mathbb{E} \bigl[ (W_{t}-W_{\tau_{t}})g(
\bar{X}_{t}) (1-\psi) \bigr]
\\
&=&\int_{\tau_t}^t\mathbb{E} \bigl[\psi
g'(\bar{X}_{t})D_u\bar{X}_{t}
\bigr] \,du+\mathbb{E} \biggl[g(\bar{X}_{t})\int_{\tau_t}^tD_u
\psi \,du \biggr]
\\
&&{}+\mathbb{E} \bigl[ (W_{t}-W_{\tau_{t}})g(
\bar{X}_{t}) (1-\psi) \bigr],
\end{eqnarray*}
where the second equality follows from the duality formula; see, for
example, Definition 1.3.1 in \cite{N}. Since for $\tau_{t}\leq u\leq t$
\begin{eqnarray*}
{\mathbb{E}} \bigl[ \psi g^{\prime}(\bar{X}_{t})D_{u}
\bar{X}_{t} \bigr] & =& {\mathbb{E}} \bigl[ \psi g^{\prime}(
\bar{X}_{t})\sigma(\bar{X}_{\tau_{t}%
}) \bigr]
\\
&=&t^{-1}
\mathbb{E} \biggl[\int_0^t D_sg(
\bar{X}_{t})\frac
{\psi
\sigma(\bar{X}_{\tau_{t}%
})}{D_s\bar{X}_t}\,ds \biggr]
\\
& =&t^{-1}{\mathbb{E}} \biggl[ g(\bar{X}_{t})\int
_{0}^{t}\psi\sigma(\bar{X}_{\tau_{t}})
\sigma^{-1} ( \bar{X}_{\tau_{s}} ) \bar{\mathcal{E}%
}_{s,t}^{-1}\delta W_{s} \biggr],
\end{eqnarray*}
one deduces
%
%
\begin{eqnarray}
\qquad\qquad\mathbb{E} [ W_{t}-W_{\tau_{t}}\rrvert\bar{X}_{t} ]
&=&t^{-1}\int_{\tau_{t}}^{t} \mathbb{E}
\biggl[ \int_{0}^{t}\psi\sigma(
\bar{X}_{\tau_{t}})\sigma^{-1} ( \bar{X}_{\tau
_{s}} )\bar{\mathcal{E}}_{s,t}^{-1}\delta W_{s}\big|
\bar{X}_{t} \biggr] \,du
\nonumber\\[-8pt]\label{espcond}\\[-8pt]
&&{}+ \mathbb{E} \biggl[\int_{\tau_t}^t
D_u\psi \,du \big| \bar{X}_{t} \biggr]+\mathbb{E} \bigl[ (
W_{t}-W_{\tau_{t}} ) (1-\psi)\rrvert\bar{X}_{t}
\bigr].\nonumber
\end{eqnarray}
Here $\delta W$ denotes the Skorohod integral. In order to obtain the
conclusion of the lemma, we need to bound the $L^p$-norm of each term
on the right-hand side of~(\ref{espcond}). In particular, we will use
the following estimate (which also proves the existence of the Skorohod
integral on the left-hand side below) which can be found in Proposition
1.5.4 in \cite{N}:
%
%
\begin{equation}
\label{controle_normep} \biggl\llVert\int_{0}^{t}
\psi\sigma(\bar{X}_{\tau_{t}})\sigma^{-1} ( \bar{X}_{\tau_{s}}
) \bar{\mathcal{E}}_{s,t}^{-1}\delta W_{s} \biggr
\rrVert_{p}\leq C(p) \bigl\llVert\psi\sigma(\bar{X}_{\tau_{t}})
\sigma^{-1} (\bar{X}_{\tau_{\cdot}} )\bar{\mathcal{E}}_{\cdot,t}^{-1}
\bigr\rrVert_{1,p},\hspace*{-35pt}
\end{equation}
where $\|F_\cdot\|_{1,p}^p =\mathbb{E} [ (\int_0^t F_s^2 \,ds )^{p/2}+
(\int_0^t \int_0^t (D_uF_s)^2 \,ds\,du )^{p/2} ]$. By Jensen's
inequality for $p\ge2$, we have
%
%
\begin{equation}
\label{upper_bound_1p} \qquad\|F_\cdot
\|_{1,p}^p \le t^{p/2-1} \int_0^t
\mathbb{E} \bigl[|F_s|^p \bigr] \,ds + t^{p-2}
\int_0^t\int_0^t
\mathbb{E} \bigl[|D_uF_s|^p \bigr]\,ds\,du
\end{equation}
and we will use this inequality to upper bound~(\ref{controle_normep}).
When $1\le p\le2$, we will use alternatively the following upper bound
$\|F_\cdot\|_{1,p}^p \le( \int_0^t \mathbb{E}[F_s^2] \,ds )^{p/2}+
(\int_0^t \int_0^t \mathbb{E}[(D_uF_s)^2] \,ds\,du )^{p/2}$ that comes
from H\"older's inequality.

For $\psi>0$, we have $\mathcal{E}_{0,t}^{-1} ( \mathcal
{E}_{0,t}%
-\bar{\mathcal{E}}_{0,t} )\leq\frac{1}{2}$ so that
$\bar{\mathcal{E}}_{0,t}\geq\frac{1}{2}\mathcal{{E}}_{0,t}>0$. From
Hypothesis~\ref{hyp_wass_pathwise}, there are constants
$0<\underline{\sigma}\le\bar{\sigma}<\infty$ such that
$0<\underline{\sigma}\leq\sigma\leq\bar{\sigma}$, and one has
\begin{eqnarray*}
&& \int_{0}^{t}{\mathbb{E}} \bigl[ \bigl( \psi
\sigma(\bar{X}_{\tau_{t}} )\sigma^{-1} ( \bar{X}_{\tau_{s}} )
\bar{\mathcal{E}}_{s,t} ^{-1} \bigr) ^{p} \bigr]\,ds
\\
&&\qquad \leq \biggl(\frac{\bar{\sigma
}}{\underline{\sigma}} \biggr)^p\int_{0}^{t}{
\mathbb{E}} \bigl[ \psi^{p}\bar{\mathcal{E}}_{0,t}^{-p}
\bar{\mathcal{E}}_{0,\eta
(s)}^{p} \bigr]\,ds
\\
&&\qquad \leq \biggl(\frac{2\bar{\sigma}}{\underline{\sigma}} \biggr)^p\sqrt
{\mathbb{E}
\bigl[\mathcal{{E}}_{0,t}^{-2p} \bigr]} \int
_0^t \sqrt{\mathbb{E} \bigl[|\bar{\mathcal{
E}}_{0,\eta
(s)}|^{2p} \bigr]}\,ds \leq C t,
\end{eqnarray*}
by using the estimates~(\ref{eq:propE}) and (\ref{eqA.13}).

Next, we focus on getting an upper bound for
%
%
\begin{equation}
\int_0^t \int_{0}^{t}{
\mathbb{E}} \bigl[ \bigl\llvert D_{u} \bigl( \psi\sigma(
\bar{X}%
_{\tau_{t}})\sigma^{-1} ( \bar{X}_{\tau_{s}}
) \bar{\mathcal{E}%
}_{s,t}^{-1} \bigr) \bigr
\rrvert^{p} \bigr] \,ds \,du.\label{eq:Dloc}%
\end{equation}
To do so, we compute the derivative using basic derivation rules, which
gives
%
%
\begin{eqnarray}\label{eq:ft}
&& D_{u} \bigl( \psi\sigma(\bar{X}_{\tau_{t}})
\sigma^{-1} ( \bar{X}%
_{\tau_{s}} ) \bar{\mathcal{
E}}_{s,t}^{-1} \bigr)\nonumber
\\
&&\qquad =D_u\psi\sigma(
\bar{X}_{\tau_{t}})\sigma^{-1} ( \bar{X}_{\tau
_{s}} )\bar{\mathcal{E}}_{s,t}^{-1}+\psi\sigma^{\prime}(
\bar{X}_{\tau_{t}})D_{u}\bar{X}_{\tau_{t}}
\sigma^{-1} ( \bar{X}_{\tau_{s}} ) \bar
{\mathcal{E}}_{s,t}^{-1}
\nonumber\\[-8pt]\\[-8pt]
&&\quad\qquad{}-\psi\sigma(\bar{X}_{\tau_{t}}) \sigma^{-2}\sigma'
( \bar{X}_{\tau_{s}} )\sigma(\bar{X}_{\tau_u}){\mathcal{
\bar E}} _{u,\tau_s}\bar{\mathcal{E}}_{s,t}^{-1}
\mathbf{1}_{u\le\tau_s}\nonumber
\\
&&\quad\qquad{}-\psi\sigma(\bar{X}_{\tau_{u}})\sigma^{-1} (
\bar{X}_{\tau_{s}} ) \bar{\mathcal{E}}_{s,t}^{-2}D_{u}
\bar{\mathcal{E}}_{s,t}.\nonumber
\end{eqnarray}
One has then to get an upper bound for the $L^p$-norm of each term. As
many of the arguments are repetitive, we show the reader only some of
the arguments that are involved. Let us start with the first term. We
have
\begin{eqnarray*}
D_u\psi&=&\varphi^{\prime} \bigl( \mathcal{E}_{0,t}^{-1}
( \mathcal{E}%
_{0,t}-\bar{\mathcal{E}}_{0,t} )
\bigr) D_{u} \bigl[ \mathcal{E}%
_{0,t}^{-1}
( \mathcal{E}_{0,t}-\bar{\mathcal{E}}_{0,t} ) \bigr]
\end{eqnarray*}
and $D_{u} [ \mathcal{E}%
_{0,t}^{-1} ( \mathcal{E}_{0,t}-\bar{\mathcal{E}}_{0,t} ) ]=
\mathcal{E}%
_{0,t}^{-2}D_{u}\mathcal{E}%
_{0,t}\bar{\mathcal{E}}_{0,t}
-\mathcal{E}%
_{0,t}^{-1}D_u\bar{\mathcal{E}}_{0,t} $. From the estimates
in~(\ref{eq:propE}), (\ref{eqA.13}) and~(\ref{eq:A13}), we obtain
%
%
\begin{equation}
\label{eq:Dest} \sup_{u\in[0,t]}\llVert D_u\psi
\rrVert_{p}\leq\bigl\|\varphi^{\prime}\bigr\|_\infty C(p).
\end{equation}
Since $\bar{\mathcal{E}}_{s,t}^{-1}=\bar{\mathcal{E}}_{0,\eta(s)}
\bar{\mathcal{E}}_{0,t}^{-1}$ and $\bar{\mathcal{E}}_{0,t}\geq
\frac{1}{2}\mathcal{{E}}_{0,t}>0$ if $\varphi^{\prime} (
\mathcal{E}_{0,t}^{-1} ( \mathcal{E}%
_{0,t}-\bar{\mathcal{E}}_{0,t} ) ) \neq0$, we have
\[
 \mathbb{E} \bigl[ \bigl\llvert D_u\psi\sigma(
\bar{X}_{\tau_{t}})\sigma^{-1} ( \bar{X}_{\tau
_{s}} )
\bar{\mathcal{E}}_{s,t}^{-1} \bigr\rrvert^p
\bigr] \le\biggl( \frac{2\bar{\sigma}}{\underline{\sigma}} \biggr)^p
\llVert
D_u\psi\rrVert_{2p}^p \mathbb{E} \bigl[
\bigl\llvert\mathcal{{E}}_{0,t}^{-1}\bar{\mathcal{
E}}_{0,\eta(s)} \bigr\rrvert^{2p} \bigr]^{1/2}.
\]
Similar bounds hold for the three other terms. Note that the highest
re\-quirements on the derivatives of~$b$ and~$\sigma$ will come from the
terms involving $D_u\bar{\mathcal{E}}$ in~(\ref{eq:ft}). Gathering all
the upper bounds,\vspace*{-2pt} we get that\break  $\llVert
\psi\sigma(\bar{X}_{\tau_{t}})\sigma^{-1} (\bar {X}_{\tau_{\cdot}}
)\bar{\mathcal{E}}_{\cdot,t}^{-1}\rrVert _{1,p}^p \le C(t^{p/2}+t^p)
\le Ct^{p/2}$ since $0\le t\le T$. From~(\ref{controle_normep}), we
finally obtain
\[
\biggl\llVert\int_{0}^{t}\psi\sigma(
\bar{X}_{\tau_{u}})\sigma^{-1} ( \bar{X}_{\tau_{s}} )
\bar{\mathcal{E}}_{s,t}^{-1}\delta W_{s} \biggr
\rrVert_{p}\leq C(p)t^{1/2}.
\]
%

We are now in position to conclude. Using Jensen's inequality, the
results (\ref{eq:propE}), (\ref{vitesse_forte}), (\ref{espcond}),
(\ref{eq:Dest}) and the definition of $\varphi$ together with
Chebyshev's inequality, we have for any $k>0$ that
\begin{eqnarray*}
&& \mathbb{E} \bigl[ \bigl\llvert\mathbb{E} \bigl[ W_{t}-W_{\tau_{t}}
|\bar{X}_{t} \bigr] \bigr|^{p} \bigr]
\\[-2pt]
&&\qquad \leq C \biggl( t^{-p}(t-\tau_{t})^{p}%
\biggl\llVert\int_{0}^{t}\psi\sigma(
\bar{X}_{\tau_{t}}%
)\sigma^{-1} ( \bar{X}_{\tau_{s}}
) \bar{\mathcal{E}}_{s,t}^{-1}\delta
W_{s} \biggr\rrVert_{p}^p
\\[-2pt]
&&\hspace*{45pt}{} +(t-\tau_{t})^{p-1}\int
_{\tau_{t}}^{t}\llVert D_u\psi\rrVert
_{p}^p\,du
\\[-2pt]
&&\hspace*{45pt}{}+\sqrt{\mathbb{E} \bigl(|W_t-W_{\tau_t}|^{2p}
\bigr)} 4^{k/2} \bigl(\mathbb{E} \bigl(|\mathcal{{E}}_{0,t}-
\bar{\mathcal{E}}_{0,t}|^{2k} \bigr)\mathbb{E} \bigl(\mathcal
{{E}}_{0,t}^{-2k} \bigr) \bigr)^{1/4} \biggr)
\\[-2pt]
&&\qquad \leq C \biggl( t^{-p/2}(t-\tau_{t})^{p}+(t-\tau
_{t})^{p}+ \biggl(\frac{1}{N} \biggr)^{ ( 2p+k ) /4}
\biggr)
\\[-2pt]
&&\qquad \leq C \biggl( \frac{1}{t^{p/2}N^{p}}+\frac{1}{N^{p/2+k/4}} \biggr).
\end{eqnarray*}\upqed
\end{pf*}\eject

\begin{pf*}{Proof of Lemma~\ref{lemme_majorations}}
The upper bounds~(\ref{eq:propE}) and (\ref{eqA.13}) on $\mathcal{{E}}$ and $\bar{\mathcal
{E}}$ are obvious since $b'$ and $\sigma'$ are bounded. Now, let us
remark that $\bar{\mathcal{E}}$ and $\mathcal{{E}}$ satisfy
\begin{eqnarray*}
\mathcal{{E}}_{u,t}&=&1+\int_u^t
\sigma' ({X}_{s})\mathcal{{E}}_{u,s}
\,dW_s+\int_u^tb'
({X}_{s})\mathcal{{E}}_{u,s} \,ds,
\\
\bar{\mathcal{E}}_{\eta_u,t}&=&1+\int_{\eta_u}^t
\sigma' (\bar{X}_{\tau_s})\bar{\mathcal{E}}_{\eta_u,\tau_s}
\,dW_s+\int_{\eta_u}^tb' (
\bar{X}_{\tau_s})\bar{\mathcal{E}}_{\eta
_u,\tau_s} \,ds.
\end{eqnarray*}
Thus, (\ref{vitesse_forte}) can be easily obtained by noticing that
$(\bar{X}_t,\bar{\mathcal{E}}_{0,t})$ is the Euler scheme for the SDE
$(X_t,\mathcal{E}_{0,t})$ which has Lipschitz coefficients, and by
using the strong convergence order of $1/2$; see, for
example,~\cite{Ka}.\vadjust{\goodbreak}

The estimate (\ref{eq:A13}) on $D_u{\mathcal{E}}$ is given, for
example, by Theorem 2.2.1 in \cite{N}. On the other hand, we have for
$\eta(s)\le u\le t$
\begin{eqnarray*}
D_u\bar{\mathcal{E}}_{\eta_s,t}&=&\sigma'(
\bar{X}_{\tau_u}) \bar{\mathcal{E}}_{\eta_s,\tau_u}
\\
&&{} +\int _{\eta_u}^t \bigl[ \sigma''(
\bar{X}_{\tau_r}) \sigma(\bar{X}_{\tau_u}) \bar{
\mathcal{E}}_{\eta_u,\tau_r} \bar{\mathcal{E}}_{\eta_s,\tau_r} +
\sigma'(\bar{X}_{\tau_r}) D_u\bar{
\mathcal{E}}_{\eta_s,\tau_r} \bigr] \,dW_r
\\
&&{}+\int_{\eta_u}^t \bigl[ b''(
\bar{X}_{\tau_r}) \sigma(\bar{X}_{\tau_u}) \bar{
\mathcal{E}}_{\eta_u,\tau_r} \bar{\mathcal{E}}_{\eta_s,\tau_r} +b'(
\bar{X}_{\tau_r}) D_u\bar{\mathcal{E}}_{\eta_s,\tau_r}
\bigr]\,dr.
\end{eqnarray*}
In order to obtain a $L^p(\Omega)$ estimate, we then
use~(\ref{eqA.13}), $b,\sigma\in C^2_b$ and Gronwall's lemma.
\end{pf*}\vspace*{-15pt}

\section{Proofs of Section~\lowercase{\protect\texorpdfstring{\ref{sec_pathwise}}{3}}}\vspace*{-5pt}\label{App_Sec2}

\begin{pf*}{Proof of Proposition \ref{prop_wass_multi}}
We use the dual representation of the Wasserstein distance
(\ref{defwas}) deduced from Kantorovitch duality theorem (see, e.g.,
Theorem 5.10, page 58 \cite{villani}),
\[
\mathcal{W}^p_p(\mu,\nu)=\sup_{\phi\in L^1(\nu)}
\biggl(\int_E\tilde{\phi}(x)\mu(dx)-\int
_E\phi(x)\nu(dx) \biggr),
\]
where $\tilde{\phi}(x)=\inf_{y\in E} (\phi(y)+|y-x|^p )$.

We also denote by $(X^{s,x}_t)_{t\in[s,T]}$ the solution to (\ref{sde})
starting from $x\in\mathbb{R}$ at time $s\in[0,T]$ and by
$(\bar{X}^{t_j,x}_t)_{t\in[t_j,T]}$ the Euler scheme starting from $x$
at time $t_j$ with $j\in\{0,\ldots,N\}$. It is enough to check that
\begin{eqnarray*}
w_k&\stackrel{\mathrm{def}}{=}&\mathcal{W}_p \bigl(\mathcal{L}
\bigl(\bar{X}_{s_1},\ldots,\bar{X}_{s_k},X^{s_k,\bar
{X}_{s_k}}_{s_{k+1}},\ldots,X^{s_k,\bar{X}_{s_k}}_{s_{n}} \bigr),
\\
&&\hspace*{21pt}
\mathcal{L} \bigl(
\bar{X}_{s_1},\ldots,\bar{X}_{s_{k-1}},X^{s_{k-1},\bar
{X}_{s_{k-1}}}_{s_{k}},\ldots,X^{s_{k-1},\bar{X}_{s_{k-1}}}_{s_{n}} \bigr) \bigr)
\end{eqnarray*}
is smaller\vspace*{1pt} than $C\sup_{0\leq t\leq
T,x\in\mathbb{R}}\mathcal{W}_p(\mathcal{L}(\bar{X}^{x}_t),\mathcal{L}(X^{x}_t))$
since $\mathcal{W}_p(\mathcal{L}(\bar{X}_{s_1},\ldots,\bar{X}_{s_n}),\allowbreak \mathcal
{L}(X_{s_1},\ldots,X_{s_n}))\leq\sum_{k=1}^nw_k$. For $f\dvtx
\mathbb{R}^n\rightarrow\mathbb{R}$ a bounded measurable function and
\[
\tilde{f}(x_1,\ldots,x_n)=\inf_{(y_1,\ldots,y_n)\in\mathbb{R}^n}
\Bigl\{ f(y_1,\ldots,y_n)+\max_{1\leq
j\leq n}|y_j-x_j|^p
\Bigr\},
\]
we set
$f_k(x_1,\ldots,x_k)=\mathbb{E}(f(x_1,\ldots,x_k,X^{s_k,x_k}_{s_{k+1}},\ldots,X^{s_k,x_k}_{s_n}))$.
First choosing
\[
(y_1,\ldots,y_{k-1},y_{k+1},\ldots,y_n)= \bigl(\bar{X}_{s_1},\ldots,\bar
{X}_{s_{k-1}},X^{s_k,y_k}_{s_{k+1}},\ldots,X^{s_k,y_k}_{s_{n}}
\bigr),
\]
then conditioning to $\sigma(W_s,s\leq s_{k})$ and using (\ref{cieds}),
next conditioning to $\sigma(W_s,s\leq s_{k-1})$ and using the dual
formulation of the Wasserstein distance, one gets
\begin{eqnarray*}
&&\mathbb{E} \bigl(\tilde{f} \bigl(\bar{X}_{s_1},\ldots,\bar
{X}_{s_k},X^{s_k,\bar
{X}_{s_k}}_{s_{k+1}},\ldots,X^{s_k,\bar{X}_{s_k}}_{s_{n}}
\bigr)
\\[-2pt]
&&\hspace*{9pt}
{}-f \bigl(\bar{X}_{s_1},\ldots,\bar{X}_{s_{k-1}},X^{s_{k-1},\bar
{X}_{s_{k-1}}}_{s_{k}},\ldots,X^{s_{k-1},\bar
{X}_{s_{k-1}}}_{s_{n}} \bigr) \bigr)
\\[-2pt]
&&\qquad \leq\mathbb{E} \Bigl(\inf_{y_k\in\mathbb{R}} \Bigl\{f \bigl(\bar
{X}_{s_1},\ldots,\bar{X}_{s_{k-1}},y_k,X^{s_k,y_k}_{s_{k+1}},\ldots,X^{s_k,y_k}_{s_{n}} \bigr)
\\[-2pt]
&&\hspace*{125pt}{}+\max_{k\leq j\leq
n}\bigl|X^{s_k,y_k}_{s_{j}}-X^{s_k,\bar{X}_{s_k}}_{s_{j}}\bigr|^p
\Bigr\}
\\[-2pt]
&&\hspace*{10pt}\quad\qquad{}-f \bigl(\bar{X}_{s_1},\ldots,\bar{X}_{s_{k-1}},X^{s_{k-1},\bar
{X}_{s_{k-1}}}_{s_{k}},\ldots,X^{s_{k-1},\bar{X}_{s_{k-1}}}_{s_{n}} \bigr) \Bigr)
\\[-2pt]
&&\qquad\leq\mathbb{E} \Bigl(\inf_{y_k\in\mathbb{R}} \bigl\{f_k(\bar
{X}_{s_1},\ldots,\bar{X}_{s_{k-1}},y_k)+C|y_k-
\bar{X}_{s_k}|^p \bigr\}
\\[-2pt]
&&\hspace*{95pt}{}
-f_k \bigl(\bar
{X}_{s_1},\ldots,\bar{X}_{s_{k-1}},X^{s_{k-1},\bar
{X}_{s_{k-1}}}_{s_{k}}
\bigr) \Bigr)
\\[-2pt]
&&\qquad\leq C\mathbb{E} \bigl(\mathcal{W}_p^p \bigl(\mathcal{L}
\bigl(X^{s_{k-1},x}_{s_k} \bigr),\mathcal{L} \bigl(\bar
{X}^{s_{k-1},x}_{s_k} \bigr) \bigr)\big|_{x=\bar{X}_{s_{k-1}}} \bigr)
\\[-2pt]
&&\qquad \leq C\sup_{x\in\mathbb{R}}\mathcal{W}^p_p \bigl(
\mathcal{L} \bigl(\bar{X}^{x}_{s_k-s_{k-1}} \bigr),\mathcal{L}
\bigl(X^{x}_{s_k-s_{k-1}} \bigr) \bigr)
\\[-2pt]
&&\qquad\leq C\sup_{0\leq t\leq T,x\in\mathbb{R}}\mathcal{W}^p_p
\bigl( \mathcal{L} \bigl(\bar{X}^{x}_t \bigr),\mathcal{L}
\bigl(X^{x}_t \bigr) \bigr).
\end{eqnarray*}\upqed
\end{pf*}

\section{Some properties of diffusion bridges}\label{diff_bridge}
Let us suppose that the SDE $dX_t=b(X_t)\,dt+\sigma(X_t)\,dW_t$,
$X_0=x$ has a~transition density $p_t(x,y)$ which is positive and of
class $\mathcal{C}^{1,2}$ with respect to $(t,x)\in\mathbb{R}_+^*
\times \mathbb{R}$. We check later in this section that this holds
under Hypothesis~\ref{hyp_wass_pathwise}. Then, the law of the
diffusion bridge with deterministic time horizon~$\mathcal{T}$ is given
by (see, e.g., Fitzsimmons, Pitman and Yor~\cite{FPY})
%
\begin{eqnarray}
\mathbb{E} \bigl[F(X_u,0\le u\le t)|X_\mathcal{T}=y \bigr]=
\mathbb{E} \biggl[F(X_u,0\le u\le t)\frac{
p_{\mathcal{T}-t}(X_t,y)}{p_\mathcal{T}(x,y)} \biggr],\nonumber
\\
\eqntext{0\le t<\mathcal{T},}
\end{eqnarray}
where $F\dvtx C([0,t],\mathbb{R})\rightarrow\mathbb{R}$ is a bounded
measurable function. Indeed for\allowbreak $g\dvtx \mathbb{R}\to\mathbb{R}$
measurable and bounded, using that $X_\mathcal{T}$ has the density
$p_\mathcal{T}(x,y)$, then the Markov property at time $t$, one checks
that
\begin{eqnarray*}
&& \mathbb{E} \biggl[\mathbb{E} \biggl[F(X_u,0\le u\le t)
\frac{
p_{\mathcal{T}-t}(X_t,y)}{p_\mathcal{T}(x,y)} \biggr] \bigg|_{y=X_\mathcal{T}
}g(X_\mathcal{T}) \biggr]
\\
&&\qquad =
\mathbb{E} \biggl[F(X_u,0\le u\le t)\int_\mathbb
{R}g(y)p_{\mathcal{T}-t}(X_t,y)\,dy \biggr]
\\
&&\qquad =\mathbb{E} \bigl[F(X_u,0\le u\le t)\mathbb{E}
\bigl[g(X_\mathcal{T})|X_t \bigr] \bigr]
\\
&&\qquad =\mathbb{E}
\bigl[F(X_u,0\le u\le t)g(X_\mathcal{T}) \bigr].
\end{eqnarray*}


We thus focus on the change of probability measure
\[
\frac{d\mathbb{P}^y}{d\mathbb{P}} \bigg|_{\mathcal{F}_t}=\frac
{p_{\mathcal{T}-t}(X_t,y)}{p_\mathcal{T}(x,y)}=:M_t,
\]
so that $\mathbb{E}[F(X_u,0\le u\le t)|X_\mathcal{T}=y]=\mathbb
{E}^y[F(X_u,0\le u\le t)]$ where $\mathbb{E}^y$ denotes the expectation
with respect to $\mathbb{P}^y$. We define $\ell_t(x,y)=\log p_t(x,y)$.
The process $(M_t)_{t\in[0,\mathcal {T})}$ is a martingale, and by
It\^o's formula, we get $dM_t=M_t \partial_x
\ell_{\mathcal{T}-t}(X_t,y) \sigma(X_t)\,dW_t$, which gives
\[
M_t=\exp\biggl( \int_0^t
\partial_x\ell_{\mathcal{T}-s}(X_s,y) \sigma
(X_s)\,dW_s -\frac{1}{2} \int
_0^t \partial_x
\ell_{\mathcal{T}
-s}(X_s,y)^2 \sigma(X_s)^2\,ds
\biggr).
\]
Girsanov's theorem then gives that for all $y\in\mathbb{R}$,
$(W^y_t=W_t-\int_0^t \partial_x\ell_{\mathcal{T}-s}(X_s,\break y)
\*\sigma(X_s)\,ds)_{t\in [0, \mathcal{T})}$ is a Brownian motion
under~$\mathbb{P}^y$, so that $(W^{X_\mathcal{T} }_t)_{t\in [0,
\mathcal{T})}$ is a Brownian motion independent of~$X_\mathcal{T}$.
Moreover, we have
%
%
\begin{equation}
\label{bridge_dyn} dX_t= \bigl[b(X_t)+
\partial_x\ell_{\mathcal{T}-t}(X_t,y)
\sigma(X_t)^2 \bigr]\,dt +\sigma(X_t)\,dW_t^y,
\end{equation}
which gives precisely the diffusion bridge dynamics.

Conversely, we would like now to reconstruct the diffusion from the
initial and the final value by using diffusion bridges. The following
result, stated in dimension one, may be generalized to higher
dimensions.

%
\begin{aprop}\label{prop_bridge} We consider an SDE
$dX_t=b(X_t)\,dt+\sigma(X_t)\,dW_t$, $X_0=x$ with a transition density
$p_t(x,y)$ positive and of class $\mathcal{C}^{1,2}$
with respect to $(t,x)\in\mathbb{R}_+^* \times\mathbb{R}$. Let $(B_t,t\ge0)$ be a
standard Brownian motion and $Z_\mathcal{T}$ be a~random variable with
density $p_\mathcal{T}(x,y)$ drawn independently from~$B$. We assume
that pathwise uniqueness holds for the SDE
%
%
\begin{eqnarray}\label{eds_pont}
dZ^{x,y}_t &=& \bigl[b
\bigl(Z^{x,y}_t \bigr)+\partial_x
\ell_{\mathcal{T}-t} \bigl(Z^{x,y}_t,y \bigr) \sigma
\bigl(Z^{x,y}_t \bigr)^2 \bigr]\,dt +\sigma
\bigl(Z^{x,y}_t \bigr)\,dB_t,\quad t\in[0,\mathcal{T}),\hspace*{-25pt}
\nonumber\\[-4pt]\\[-4pt]
Z^{x,y}_0 &=& x,\nonumber\hspace*{-25pt}
\end{eqnarray}
for any $x,y\in\mathbb{R}$, and set $Z_t= Z^{x,Z_\mathcal{T}}_t$ for $t
\in
[0,\mathcal{T})$. Then, $(Z_t)_{t \in[0,\mathcal{T}]}$ and
$(X_t)_{t\in[0,\mathcal{T}]}$ have the same law.
\end{aprop}

A consequence of this result is that $(Z_t,t \in[0,\mathcal{T}])$ has
continuous paths, which gives that $\lim_{t\rightarrow
\mathcal{T}-}Z^{x,y}_t = y$ a.s., $dy$-a.e.

\begin{pf}
Let $t\in[0,\mathcal{T})$ and $F\dvtx C([0,t],\mathbb
{R})\to\mathbb{R}$ and $g:\mathbb{R}\to\mathbb{R}$ be bounded and
measurable functions. Since pathwise uniqueness for the
SDE~(\ref{eds_pont}) implies weak uniqueness, we get
\begin{eqnarray*}
\mathbb{E} \bigl[F \bigl(Z_{u}^{x,y},0\leq u\leq t \bigr)
\bigr]&=&\mathbb{E}^y \bigl[F(X_u,0\leq u\leq t) \bigr]
\\
&=& \mathbb{E} \biggl[F(X_u,0\leq u\leq t)\frac
{p_{\mathcal{T}
-t}(X_t,y)}{p_\mathcal{T}(x,y)} \biggr].
\end{eqnarray*}
Thus we have
\begin{eqnarray*}
\mathbb{E} \bigl[F(Z_{u},0\leq u\leq t)g(Z_\mathcal{T})
\bigr]&=&\mathbb{E} \biggl[F(X_u,0\leq u\leq t)\int
_\mathbb{R}p_{\mathcal{T}
-t}(X_t,y)g(y)\,dy \biggr]
\\
&=&
\mathbb{E} \bigl[F(X_u,0\leq u\leq t)g(X_\mathcal{T} )
\bigr].
\end{eqnarray*}
Hence the finite-dimensional marginals of the two processes are equal.
Since $(X_t)_{t\in[0,\mathcal{T}]}$ has continuous paths and
$(Z_t)_{t\in[0,\mathcal{T}]}$ has c\`adl\`ag paths (continuous on
$[0,\mathcal{T})$ with a possible jump at $\mathcal{T}$), this
completes the proof.
\end{pf}


From now on, we assume that Hypothesis~\ref{hyp_wass_pathwise} holds.
We introduce the Lamperti transformation of the stochastic process
$(X_t,t\ge0)$. We\vspace*{-4pt} define $\varphi(x)=\int_0^x\frac{dy}{\sigma(y)}$ and
$\alpha(y)= (\frac{b}{\sigma}-\frac{\sigma'}{2} )\circ
\varphi^{-1}(y)$, $\hat{X}_t\stackrel{\mathrm{def}}{=}\varphi(X_t)$ so
that we have
%
%
\begin{equation}
\label{sde_lamperti} d\hat{X}_t= \alpha(
\hat{X}_t)\,dt + dW_t,\qquad t\in[0,T].
\end{equation}
By\vspace*{-2pt} Hypothesis~\ref{hyp_wass_pathwise}, $\varphi$ is a $C^5$ bijection,
$\alpha\in C^3_b$ and both $\varphi$ and $\varphi^{-1}$ are Lipschitz
continuous. We denote by $\hat{p}_t(\hat{x},\hat{y})$ the transition density
of~$\hat{X}$ and $\hat{\ell}_t(\hat{x},\hat{y})=\log(\hat
{p}_t(\hat{x},\hat{y}))$.


%
\begin{alem}\label{lem_densite} The density $\hat{p}_t(\hat{x},\hat
{y})$ is
$C^{1,2}$ with respect to $(t,\hat{x})\in\mathbb{R}_+^*\times
\mathbb{R}$. Besides, we have
\[
\partial_{\hat{x}} \hat{\ell}_t(\hat{x},\hat{y})=
\frac{\hat
{y}-\hat{x}}{t}-\alpha(\hat{x})+ g_t(\hat{x},\hat{y}),
\]
where $g_t(\hat{x},\hat{y})$ is a continuous function on $\mathbb
{R}_+\times
\mathbb{R}^2$ such that
$\partial_{\hat{x}} g_t(\hat{x},\hat{y})$ and $\partial_{\hat{y}}
g_t(\hat{x},\hat{y})$
exist and
\[
\forall T>0\qquad
\sup_{t\in[0,T], \hat{x},\hat{y}\in\mathbb
{R}}\bigl|\partial_{\hat{x}}
g_t(\hat{x},\hat{y})\bigr|+\bigl|\partial_{\hat{y}} g_t(\hat
{x},\hat{y})\bigr|<\infty.
\]
\end{alem}

\begin{pf}
It is well known that we can express the transition
density~$\hat{p}_t(\hat{x},\hat{y})$ by using Girsanov's theorem as an
expectation on a Brownian bridge between $\hat{x}$~and~$\hat{y}$.
Namely, since $\alpha$ and its derivatives are bounded, we can apply a
result stated in Gihman and Skorohod~\cite{gs} (Theorem~1, Chapter~3,
Section~13) or in Rogers \cite{rog} to get that $\hat{p}_t(\hat{x},\hat{y})$ is positive and
\begin{eqnarray*}
\hat{\ell}_t(\hat{x},\hat{y})&=&-\frac{(\hat{x}-\hat
{y})^2}{2t}+\int
_{\hat{x}}^{\hat{y}
}\alpha(z)\,dz
\\
&&{} +\log\mathbb{E}
\biggl(\exp\biggl({-\frac{1}{2}\int_0^t\bigl(\alpha
'+\alpha
^2\bigr)\biggl(\hat{x}+W_s+\frac{s}{t}(\hat{y}-\hat{x}-W_t)\biggr)\,ds} \biggr)\biggr)
\\
&&{} -\frac{1}{2}\log(2\pi t).
\end{eqnarray*}
Clearly, $\hat{\ell}_t(\hat{x},\hat{y})$ is $C^{1,2}$ in
$(t,\hat{x})\in\mathbb{R}_+^*\times
\mathbb{R}$ (we can use carefree the dominated convergence theorem for
the third
term since $\alpha\in C^3_b$), and we have
\begin{eqnarray*}
g_t(\hat{x},\hat{y})
&=& -\frac{1}{2}
\biggl(\mathbb{E} \biggl[\exp\biggl({-\frac
{1}{2}\int_0^t\bigl(\alpha'+\alpha^2\bigr)\biggl(\hat{x}+W_s+\frac{s}{t}(\hat{y}-\hat
{x}-W_t)\biggr)\,ds}\biggr)
\\
&&\hspace*{33pt}{}\times \int_0^t\frac{t-s}{t}\bigl(\alpha''+2\alpha\alpha'\bigr)\biggl(\hat
{x}+W_s+\frac
{s}{t}(\hat{y}-\hat{x}-W_t)\biggr)\,ds \biggr]\biggr)
\\
&&\hspace*{9pt}{}\bigg/\biggl({\mathbb{E} \biggl[\exp\biggl({-\frac
{1}{2}\int
_0^t\bigl(\alpha'+\alpha^2\bigr)\biggl(\hat{x}+W_s+\frac{s}{t}(\hat{y}-\hat
{x}-W_t)\biggr)\,ds}\biggr) \biggr]}\biggr).
\end{eqnarray*}
This is a continuous function on~$\mathbb{R}_+\times\mathbb{R}^2$,
and we easily
conclude by using the dominated convergence theorem and~$\alpha\in
C^3_b$.
\end{pf}

By straightforward calculations, we have
\[
p_t(x,y)=\frac{1}{\sigma(y)}\hat{p}_t \bigl(\varphi(x),
\varphi(y) \bigr)
\]
and $p_t(x,y)$ is thus positive and $C^{1,2}$ with respect to $(t,x)$.
The diffusion bridge~(\ref{bridge_dyn}) is thus well defined. Since
$\partial_x \ell_t(x,y) =\frac{1}{\sigma(x)}\partial_{\hat x}
\hat{\ell}_t(\varphi(x),\varphi(y))$, we get by It\^{o} formula
from~(\ref{bridge_dyn})
\begin{eqnarray*}
d \hat{X}_t&=& \bigl[\alpha(\hat{X}_t)+
\partial_{\hat{x}} \hat{\ell}_{\mathcal{T}-t} \bigl(\hat{X}_t,
\varphi(y) \bigr) \bigr]\,dt+dW^y_t,
\\
dW^y_t &=&dW_t-
\partial_{\hat{x}} \hat{\ell}_{\mathcal{T}-t} \bigl(\hat{X}_t,
\varphi(y) \bigr)\,dt.
\end{eqnarray*}
Therefore, as one could expect, the Lamperti transform on the diffusion
bridge coincides with the diffusion bridge on the Lamperti transform.

%
\begin{aprop}\label{prop_bridge2}
Let Hypothesis~\ref{hyp_wass_pathwise} hold. There exists a
deterministic constant~$C$ such that
\[
\forall\mathcal{T}\in(0,T], x,x',y,y'\in\mathbb{R}\qquad
\sup_{t\in[0,\mathcal{T})} \bigl|Z^{x,y}_t-Z^{x',y'}_t\bigr|
\le C \bigl(\bigl|x-x'\bigr|\vee\bigl|y-y'\bigr| \bigr)
\]
and in particular, pathwise uniqueness holds for~(\ref{eds_pont}).
\end{aprop}

\begin{pf}
For $\hat{x},\hat{y}\in\mathbb{R}$, we consider the following SDE:
%
%
\begin{eqnarray}\label{pont_Z2}
d \hat{Z}^{\hat{x},\hat
{y}}_t&=&dB_t+
\biggl[\frac{\hat{y}-\hat{Z}^{\hat{x},\hat{y}}_t}{\mathcal
{T}-t}+g_{\mathcal{T}-t} \bigl(\hat{Z}^{\hat{x},\hat{y}}_t,
\hat{y} \bigr) \biggr]\,dt,\qquad t\in[0, \mathcal{T}),
\nonumber\\[-20pt]\\
\hat{Z}^{\hat{x},\hat{y}}_{0}&=&\hat{x},\nonumber
\end{eqnarray}
which corresponds to the diffusion bridge on the Lamperti
transform~$\hat{X}$. We set
$\Delta_t=\hat{Z}^{\hat{x},\hat{y}}_t-\hat{Z}^{\hat{x}',\hat
{y}'}_t$ for
$t\in[0,\mathcal{T})$ and $\hat{x}',\hat{y}' \in\mathbb{R}$. We have
\[
d \Delta_t = \biggl[\frac{\hat{y}-\hat{y}'-\Delta_t}{\mathcal{T}-t}
+g_{\mathcal{T}-t} \bigl(
\hat{Z}^{\hat{x},\hat{y}}_t,\hat{y} \bigr)-g_{\mathcal{T}-t} \bigl(\hat
{Z}^{\hat{x}',\hat{y}'}_t,\hat{y}' \bigr) \biggr]\,dt
\]
and thus $d(|\Delta_t| \vee|\hat{y}-\hat{y}'|)=
\operatorname{sign}(\Delta_t)\mathbf{1}_{|\Delta_t|\ge|\hat{y}-\hat{y}'|}\,d\Delta_t$.
On the one hand, we observe that
$\mathbf{1}_{|\Delta_t|\ge|\hat{y}-\hat{y}'|}[\operatorname{sign}(\Delta_t)
(\hat{y}-\hat{y}')-|\Delta_t|]\le0$. On the other hand, $g_t$ is
uniformly Lipschitz w.r.t. $(\hat{x},\hat{y})$ on $t\in[0,T]$ by
Lemma~\ref{lem_densite}, which leads to
\[
d \bigl(|\Delta_t| \vee\bigl|\hat{y}-\hat{y}'\bigr| \bigr)\le C
\bigl(|\Delta_t| \vee\bigl|\hat{y}-\hat{y}'\bigr| \bigr)
\]
for some positive constant~$C$. Gronwall's lemma gives then
$|\Delta_t|\le e^{CT}(|\hat{x}-\hat{x}'|\vee|\hat{y}-\hat{y}'|)$.
This gives in
particular pathwise uniqueness for~(\ref{pont_Z2}).

Now, let us\vspace*{-1pt} assume that $(Z^{x,y}_t)_{t\in[0,\mathcal{T})}$
solves~(\ref{eds_pont}). Then
$\varphi(Z^{x,y}_t)$ solves~(\ref{pont_Z2}) with $\hat{x}=\varphi
(x)$ and
$\hat{y}=\varphi(y)$, and we necessarily have
$Z^{x,y}_t=\varphi^{-1}(\hat{Z}_t^{\varphi(x),\varphi(y)})$ by pathwise
uniqueness. Both $\varphi$ and $\varphi^{-1}$ are Lipschitz, and we
denote by $K$ a~common Lipschitz constant. Then we get
\begin{eqnarray*}
\bigl|Z^{x,y}_t-Z^{x',y'}_t\bigr| &=&\bigl|
\varphi^{-1} \bigl(\hat{Z}_t^{\varphi(x),\varphi(y)} \bigr) -
\varphi^{-1} \bigl(\hat{Z}_t^{\varphi(x'),\varphi(y')} \bigr) \bigr|
\\
&\le&
K^2 e^{CT} \bigl(\bigl|x-x'\bigr|\vee\bigl|y-y'\bigr|
\bigr),
\end{eqnarray*}
which gives the desired result.
\end{pf}
\end{appendix}




\printaddresses

\end{document}